\documentclass[12pt,reqno]{amsart} 
\usepackage{amssymb}\usepackage{epsf}

\numberwithin{equation}{section}

\theoremstyle{definition}
\DeclareMathOperator*{\adj}{adj}
\DeclareMathOperator*{\ind}{ind}

\def\stackreb#1#2{\ \mathrel{\mathop{#1}\limits_{#2}}}
\DeclareMathOperator{\Li}{Li}
\newcommand{\eg}{\Gamma}
\newcommand{\beq}{\begin{equation}}
\newcommand{\be}{\begin{equation}}
\newcommand{\ee}{\end{equation}}
\newcommand{\ba}{\begin{eqnarray}}
\newcommand{\ea}{\end{eqnarray}}

\newcommand{\lab}[1]{\label{#1}}

\newcommand{\C}{\mathbb C}
\newcommand{\R}{\mathbb R}
\newcommand{\Z}{\mathbb Z}

\newcommand{\T}{\mathbb T}

\newcommand{\ve}{\varepsilon}

\begin{document}

\title
{Introduction to the theory of elliptic  hypergeometric integrals
}

\author{Vyacheslav P. Spiridonov}

\address{Laboratory of theoretical physics, Joint Institute for Nuclear
Research, Dubna, Moscow reg., 141980, Russia and
National Research University Higher School of Economics, Moscow, Russia}

\begin{abstract}
We give a brief account of the key properties of elliptic hypergeometric integrals --- a relatively
recently discovered top class of transcendental special functions of hypergeometric type.
In particular, we describe an elliptic generalization of Euler's and Selberg's beta integrals,
elliptic analogue of the Euler-Gauss hypergeometric function and some multivariable elliptic
hypergeometric functions on root systems. The elliptic Fourier transformation and
corresponding integral Bailey lemma technique is outlined together with a
connection to the star-triangle relation and Coxeter relations for a permutation group.
We review also the interpretation of elliptic hypergeometric integrals as superconformal indices
of four dimensional supersymmetric quantum field theories and corresponding applications
to Seiberg type dualities.
\end{abstract}

\maketitle

\tableofcontents

\section{Introduction}

The Euler-Gauss hypergeometric function \cite{aar} is one of the most useful
classical special functions. Its most popular definition is given by the ${}_2F_1$-series:
\beq
F(a,b;c;x):={}_2F_1(a,b;c;x)
=\sum_{n=0}^\infty \frac{(a)_n(b)_n}{n!(c)_n}x^n,
\qquad |x|<1,
\lab{2F1}\ee
where $(a)_n=a(a+1)\ldots (a+n-1)$ is the Pochhammer symbol.
Alternatively, it can be defined by the Euler integral representation
\beq
F(a,b;c;x)=\frac{\Gamma(c)}{\Gamma(c-b)\Gamma(b)}\int_0^1 t^{b-1}(1-t)^{c-b-1}
(1-xt)^{-a}dt,
\lab{eulerint}\ee
where $\text{Re}(c)>\text{Re}(b)>0$ and $x\notin [1,\infty[ $,
or the Barnes integral representation
\beq
F(a,b;c;x)=\frac{\Gamma(c)}{\Gamma(a)\Gamma(b)}\int_{-i\infty}^{i\infty}
\frac{\Gamma(a+u)\Gamma(b+u)\Gamma(-u)}{\Gamma(c+u)}(-x)^u du,
\lab{barnesint}\ee
where  $\Gamma(x)$ is the Euler gamma function
$$
\Gamma(x)=\int_0^\infty t^{x-1}e^{-t}dt, \qquad \text{Re}(x)>0.
$$
In \eqref{barnesint} the poles of the integrand $u=-a-k, -b-k$ and $u=k$, $k\in\Z_{\geq0}$, are separated
by the integration contour.

The function $F(a,b;c;x)$ satisfies a special differential equation called the hypergeometric equation:
\beq
x(1-x)y''(x)+(c-(a+b+1)x)y'(x)-ab y(x)=0,
\lab{HE}\ee
determining the solution analytic around the regular singular point $x=0$. This is the second
order differential equation with three regular singularities fixed at $x=0, 1, \infty$
by linear fractional transformation.

For $x=1$ the value of function $F(a,b;c;1)$ \eqref{eulerint} can be computed due to
the explicit evaluation of Euler's beta integral:
\beq
\int_0^1t^{x-1}(1-t)^{y-1}dt=\frac{\Gamma(x)\Gamma(y)}
{\Gamma(x+y)}, \qquad \text{Re}(x), \text{Re}(y)>0.
\label{beta}\ee

Restrictions on the parameters indicated above lead to well defined functions,
they may be relaxed by analytic continuation.

All these exact formulas and related ones
were generalized in many different ways. We mention the most essential
developments:
\begin{itemize}

\item  extension to higher order hypergeometric functions $_{n+1}F_n$,

\item $q$-deformation of plain hypergeometric functions,

\item  extension of univariate to multivariable special functions,

\item elliptic deformation of all above functions.

\end{itemize}

It is the last step which will be our main subject in these notes.
It represents a relatively recent development in the theory of
special functions with the basic results obtained around 2000.
For describing the most general elliptic hypergeometric functions
one has to use integral representations \cite{spi:umn}, since the infinite series
of the corresponding type are not well defined.
Note however, that the first examples of elliptic
hypergeometric functions emerged in the terminating series
form as particular elliptic function solutions of the Yang-Baxter equation \cite{FT}
which were constructed in a case-by-case manner in \cite{djkmo}.

The most interesting elliptic hypergeometric integrals are
associated with two independent root systems related in a remarkable way
to supersymmetric quantum field theories, where these integrals emerge as
superconformal indices \cite{RR2016}. The first root system determines their structure
as matrix integrals over the Haar measure of a particular compact Lie group
(the gauge group in field theory), and the second one is related to a Lie group
of symmetry transformations of functions in parameters (the flavor group in field theory).
There are many exact relations between
such integrals, a large number of which are still in a conjectural form.

We shall not try to cover all aspects of the theory, but consider some
introductory material at the elementary level and give a brief review
of more recent developments. There are other surveys on this subject
\cite{rosengren}, \cite{ros-war}, \cite{spi:umnrev}, where some of
the skipped topics are discussed. A deep algebraic geometry point
of view on the functions of interest is given in \cite{rai:noncom}.

\section{Elliptic hypergeometric integrals}

The very first basic example of elliptic hypergeometric integrals was discovered in \cite{spi:umn}.
Let us start from the conceptual definition of  such integrals
introduced in \cite{spi:theta2}. For simplicity we limit its consideration only to the univariate case.

The key property of the univariate elliptic hypergeometric integrals
is that they are defined as contour integrals
$$
I:=\int_C\Delta(u)du,
$$
whose kernel $\Delta(u)$ satisfies a first order finite difference equation
\begin{equation}
\Delta(u+\omega_1)=f(u;\omega_2,\omega_3)\Delta(u),
\label{defEHI}\end{equation}
where the coefficient $f(u;\omega_2,\omega_3)$
is an elliptic function with periods $\omega_2$ and $\omega_3$,
and $\omega_{1,2,3}$ are some incommensurate complex numbers.
Incommensurability means that $\sum_{k=1}^3n_k\omega_k\neq 0$ for
$n_k\in\mathbb{Z}$.

Elliptic functions form a particular beautiful family of special functions \cite{akh}.
Let us remind that they are defined as the meromorphic doubly-periodic functions:
$$
f(u+\omega_2)=f(u+\omega_3)=f(u), \quad  \text{Im}(\omega_2/\omega_3)\neq 0.
$$
Consider their general structure before discussing solutions of the
defining equation \eqref{defEHI}. For that we need an infinite product
$$
(z;p)_\infty:=\prod_{j=0}^\infty(1-zp^j),\quad |p|<1, \, z\in\mathbb{C}.
$$
With its help we define a Jacobi theta function as
$$
\theta(z;p):=(z;p)_\infty(pz^{-1};p)_\infty,\quad
z\in\ \mathbb{C}^\times.
$$
It has important symmetry properties:
\beq
\theta(pz;p)=\theta(z^{-1};p)=-z^{-1}\theta(z;p).
\lab{thetasym}\ee
Using the Jacobi triple product identity one can write the Laurent series expansion
\beq
\theta(z;p)=\frac{1}{(p;p)_\infty} \sum_{k\in\Z} (-1)^kp^{k(k-1)/2}z^k.
\lab{Jtri}\ee
For convenience we provide the standard odd Jacobi theta function definition:
\begin{eqnarray} \lab{theta1} &&
 \theta_1(u|\tau)=-\theta_{11}(u)=-\sum_{k\in \mathbb Z}e^{\pi {i} \tau(k+1/2)^2}
e^{2\pi {i} (k+1/2)(u+1/2)}
\\  && \makebox[4em]{}
={i}p^{1/8}e^{-\pi {i} u}(p;p)_\infty\theta(e^{2\pi {i} u};p),
\quad p=e^{2\pi {i} \tau}.
\nonumber\end{eqnarray}

For interested readers  we suggest small calculational tasks like the following one. {\em Exercise:}
find zeros of $\theta(z;p)$,  verify \eqref{thetasym},
deduce the general quasiperiodicity relation for $\theta(p^kz;p),\, k\in\Z$,
and prove identity \eqref{Jtri}. The latter problem can be solved
by computing the $z$-series expansion for the symmetric finite product
$\prod_{k=1}^n(1-zp^{k-1/2})(1-z^{-1}p^{k-1/2})$
(using the $z\to pz$ functional equation for that) and  taking the limit $n\to\infty$.

There is nice factorized representation of the elliptic functions in terms
of the Jacobi theta functions. Denote
$$
p := e^{2 \pi {i} \omega_3/\omega_2}, \quad  z:=e^{2\pi{i} u/\omega_2}.
$$
Then, according to the theorem established by Abel and Jacobi, one can write
up to a multiplicative constant,
\beq
h(z;p):=f(u;\omega_2,\omega_3)=\prod_{k=1}^m\frac{\theta(t_kz;p)}{\theta(w_kz;p)},
\qquad \prod_{k=1}^mt_k=\prod_{k=1}^mw_k.
\label{elliptic}\ee
Indeed, the periodicity $f(u+\omega_2)=f(u)$
of this meromorphic function is evident. Since the shift
$u\to u+\omega_3$ is equivalent to $z\to pz$, one has
$f(u+\omega_3)/f(u)=h(pz)/h(z)=\prod_{k=1}^m(w_k/t_k)=1$.
So, we have an elliptic function with $m$ poles (or zeros) in the
fundamental parallelogram of periods $(\omega_2,\omega_3)$.
Vice versa, given an elliptic function with $m$
poles and zeros at the points fixed by parameters $t_k$ and $w_k$, we can divide
it by \eqref{elliptic} and see that the resulting function is bounded on
$\mathbb{C}$ (it is doubly periodic and has no poles), and by the Liouville theorem
it is constant. The parameter $m$ is called the order
of elliptic functions and we call the linear constraint
$\prod_{k=1}^mt_k=\prod_{k=1}^mw_k$
the balancing condition (it explains the origin of the old
notion of balancing in the theory of hypergeometric functions \cite{gas-rah:basic,spi:theta1}).

It is convenient to use compact notation
$$
\theta(a_1,\ldots,a_k;p):=\theta(a_1;p)\cdots\theta(a_k;p), \quad
\theta(at^{\pm 1};p):=\theta(at;p)\theta(at^{-1};p).
$$
Then the ``addition'' formula for theta functions takes the form
\begin{eqnarray}
\theta(xw^{\pm 1},yz^{\pm 1};p) -\theta(xz^{\pm 1},yw^{\pm 1};p)
=yw^{-1}\theta(xy^{\pm 1},wz^{\pm 1};p).
\label{addition}\end{eqnarray}
The proof of this relation is rather easy. The ratio of the left-
and right-hand sides satisfies equation $h(px)=h(x)$ (i.e., it is
$p$-elliptic) and represents a bounded function of the variable $x\in\C^\times $.
Therefore it does not depend on $x$ according to the Liouville theorem,
but for $x=w$ the equality is evident.

Let us turn now to the elliptic hypergeometric integrals.
In terms of the multiplicative coordinate $z=e^{2\pi{i} u/\omega_2}$
elliptic functions are determined by the equation $h(pz)=h(z)$. Now we demand that
the integrand $\Delta(u)=:\rho(z)$ is a meromorphic function of
$z\in\mathbb{C}^\times$, which is an additional strong restriction. Then it is convenient to
introduce a second base variable $q:=  e^{2 \pi {i}\omega_1/\omega_2},$
so that the shift $u\to u+\omega_1$ becomes equivalent to the
multiplication $z\to qz$. Changing the integration
variable, we come to the following definition of elliptic hypergeometric integrals:
\begin{equation}
I_{EHI}=\int\rho(z)\frac{dz}{z}, \qquad \rho(qz)=h(z;p)\rho(z),\quad h(pz)=h(z),
\label{EHImult}\end{equation}
where the explicit form of $h(z;p)$ is given in \eqref{elliptic}.

Because of the factorized form of $h(z;p)$, for  solving
the equation $\rho(qz)=h(z;p)\rho(z)$ it is sufficient to
solve the linear first order $q$-difference equation
with a simple theta function coefficient
\beq
\gamma(qz)=\theta(z;p)\gamma(z).
\lab{ellgammaeq}\ee
One can check that a particular solution of \eqref{ellgammaeq} is given by
the function
\beq
\gamma(z)=\Gamma(z;p,q)
:=\prod_{j,k=0}^\infty\frac{1-z^{-1}p^{j+1}q^{k+1}}{1-zp^{j}q^{k}},\qquad |p|, |q|<1,
\lab{eg}\ee
which is called the (standard) elliptic gamma function.
Note that the original equation \eqref{ellgammaeq} does not impose the constraint $|q|<1$,
whereas in \eqref{eg} we have such an additional restriction.

The problem of generalizing the Euler gamma function was considered
by Barnes, who defined the multiple gamma functions of arbitrary order \cite{bar:multiple},
and Jackson \cite{jac:basic}, who introduced the basic versions of gamma functions.
Although the function \eqref{eg} is related to their considerations, its usefullness
was established only in modern time after the work of Ruijsenaars \cite{rui:first},
where the term ``elliptic gamma function'' was introduced.
A further systematic investigation of this function
was performed by Felder and Varchenko \cite{fel-var:elliptic} who discovered
its $\textrm{SL}(3,\mathbb{Z})$ symmetry transformations (they also
pointed out that this function appeared implicitly already in Baxter's work
on the eight vertex model \cite{bax}). In \cite{spi:theta2} the author
constructed the modified elliptic gamma function, which gives
a solution of equation \eqref{ellgammaeq} in the regime $|q|=1$
(it is meromorphic in $u\propto \log z$, not $z$). It will be described in the
next section.

{\em Exercise}: derive the solution \eqref{eg}
from scratch by iterations using the factorized form of $\theta(z;p)$.

Changing the variable $z\to t_kz, w_kz$ in \eqref{ellgammaeq}, we
find solutions of the equation defining $\rho(z)$ for each theta function
factor in $h(z;p)$. The final result is evident now: the general
univariate elliptic hypergeometric integral has the form
\beq
I_{EHI}(\underline{t},\underline{w};p,q)=\int\prod_{k=1}^{m}
\frac{\Gamma(t_kz;p,q)}{\Gamma(w_kz;p,q)}\frac{dz}{z},
\qquad \prod_{k=1}^mt_k=\prod_{k=1}^mw_k,
\lab{EHI}\ee
where one has to specify the contour of integration. The typical
choice is a closed contour encircling the essential singularity
point $z=0$, e.g. the unit circle $\mathbb{T}$.
Surprisingly, these functions generalize all previously known univariate
ordinary and $q$-hypergeometric functions. They depend on $2m$ complex
variables $t_k$ and $w_k$ subject to one constraint. We shall not describe
explicitly the limits to lower level hypergeometric objects, but
only indicate how it can be done. Note that for $0<|p|, |q|<1$ it is not
possible to simplify functions \eqref{EHI} by taking parameters to zero
or infinity. Therefore all well defined degenerations require limits to
the boundary values of bases.

{\em Exercise:} investigate the $p\to 0$ limit of \eqref{EHI} for fixed
parameters and when some of the parameters behave as powers of $p$.

Consider the uniqueness of the derived expression for $I_{EHI}$.
Evidently, solutions of equation \eqref{ellgammaeq} are defined up to
the multiplication by an arbitrary elliptic function of some order $l$,
whose general form was fixed in \eqref{elliptic}. However, one can write
$$
\prod_{j=1}^l\frac{\theta(a_jz;p)}{\theta(b_jz;p)}
=\prod_{j=1}^l\frac{\Gamma(qa_jz;p,q)\Gamma(b_jz;p,q)}
{\Gamma(a_jz;p,q)\Gamma(qb_jz;p,q)}, \qquad \prod_{j=1}^la_j=\prod_{j=1}^lb_j,
$$
and see that the right-hand side expression can be absorbed to the
kernel in \eqref{EHI} by extension of the set of parameters $\{t_k\}\to
\{t_k, qa_j,b_j\}$ and $\{w_k\}\to\{w_k, a_j,qb_j\}$ without violating
the balancing condition. Therefore \eqref{EHI} can be considered as
a general solution. As to the initial equation \eqref{defEHI},
its solutions can be multiplied by arbitrary function of period
$\omega_1$ which cannot be fixed without imposing additional constraints.

\section{Properties of the elliptic gamma function}

For describing properties of the elliptic gamma function
we take the same $\omega_{1,2,3}$ as in the previous section
and introduce three  bases
$$
p \ = \ e^{2 \pi {i}
\omega_3/\omega_2}, \quad q \ = \ e^{2 \pi {i}
\omega_1/\omega_2}, \quad r \ = \ e^{2 \pi {i} \omega_3/\omega_1}
$$
and their particular modular partners
$$
\widetilde{p} \ = \ e^{-2 \pi {i}
\omega_2/\omega_3}, \quad \widetilde{q} \ = \ e^{-2 \pi {i}
\omega_2/\omega_1}, \quad \widetilde{r} \ = \ e^{-2 \pi {i}
\omega_1/\omega_3}.
$$

The first relation we draw attention to is an evident symmetry in bases
$$
\Gamma(z;p,q)=\Gamma(z;q,p),
$$
which looks quite unexpected taking into account how asymmetrically
the bases $p$ and $q$ enter equation \eqref{ellgammaeq}.
Due to this symmetry one actually has two finite-difference equations
$$ 
\Gamma(qz;p,q)=\theta(z;p)\Gamma(z;p,q),
\quad\Gamma(pz;p,q)=\theta(z;q)\Gamma(z;p,q).
$$

Poles and zeros of the elliptic gamma function form a two-dimensional array of
geometric progressions
$$
z_{\mathrm{poles}}=p^{-j}q^{-k},\quad
z_{\mathrm{zeros}}=p^{j+1}q^{k+1},\quad j,k\in\mathbb{Z}_{\geq 0}.
$$
The inversion relation for $\Gamma(z;p,q)$ has the form
\begin{equation}
\Gamma(z;p,q)=\frac{1}{\Gamma(\frac{pq}{z};p,q)},
\label{inv}\end{equation}
and there is a useful normalization condition $\Gamma(\sqrt{pq};p,q)=1.$

The quadratic transformation
$$
\Gamma(z^2;p,q)=\Gamma(\pm z,\pm q^{1/2}z,\pm p^{1/2}z,\pm (pq)^{1/2}z;p,q)
$$
can be established by a direct analysis of the infinite products.
Here and below we use the conventions
\begin{eqnarray*}
&& \Gamma(t_1,\ldots,t_k;p,q):=\Gamma(t_1;p,q)\cdots\Gamma(t_k;p,q),\quad
\\ &&
\Gamma(\pm z;p,q):=\Gamma(z;p,q)\Gamma(-z;p,q),
\\ &&
\Gamma(tz^{\pm k};p,q):=\Gamma(tz^k;p,q)\Gamma(tz^{-k};p,q).
\end{eqnarray*}
The limiting relation
\beq
\lim_{z\to1}(1-z)\Gamma(z;p,q)=\frac{1}{(p;p)_\infty (q;q)_\infty}
\lab{limit}\ee
is required for residue calculus and reduction of integrals
to terminating elliptic hypergeometric series (non-terminating such
series do not converge).

Taking the logarithm of the infinite product \eqref{eg}, expanding the
logarithms of individual factors $\log(1-x)=-\sum_{n=1}^\infty x^n/n,$
and changing the summation order yields the following representation
\beq
\Gamma(z;p,q)=\exp\left(\sum_{n=1}^\infty \frac{1}{n}
\frac{z^n-(pq/z)^n}{(1-p^n)(1-q^n)}\right),
\lab{egexp}\ee
which converges for $|pq|<|z|<1$ and is very useful for quantum field theory purposes.

Denote $p=e^{-\delta}$ and consider the limit $\delta \to 0$.
The leading asymptotics takes the form
\begin{equation}
\Gamma(z;p,q)=\exp\left(\frac{1}{\delta}E_2(z;q)-\frac{1}{2}\log\theta(z;q)\right)(1+O(\delta)),
\label{asymp}\end{equation}
where
$$
E_2(z;q)=\sum_{n=0}^\infty \Li_2(q^nz)- \sum_{n=1}^\infty \Li_2(q^n/z),
\quad \Li_2(z)=\sum_{n=1}^\infty \frac{z^n}{n^2}.
$$
$\Li_2(z)$ is known as Euler's dilogarithm function and $E_2(z;q)$ is
directly related to the  elliptic dilogarithm function, which
recently emerged in the computation of a sunset Feynman diagram \cite{BKV}
as the difference $\hat E(z;q)=E(z;q)- E(-z;q)$. The latter function emerges in
the asymptotics of the ratio $\Gamma(z;p,q)/\Gamma(-z;p,q)$. A different relation between
the elliptic gamma function with $p=q$ and the elliptic dilogarithm was described in \cite{PZ}.

We shall need also the second order generalization of the elliptic gamma function
$$
\Gamma(z;p,q,t)=\prod_{j,k,l=0}^\infty (1-zp^jq^kt^l)
(1-z^{-1} p^{j+1}q^{k+1}t^{l+1}), \quad |t|,|p|,|q|<1,\; z\in\mathbb{C}^\times.
$$
It is related to function \eqref{eg} via the difference equation
\beq
\Gamma(qz;p,q,t)=\Gamma(z;p,t)\Gamma(z;p,q,t),
\label{2ndordereqn}\ee
and its inversion relation has the form $\Gamma(pqtz;p,q,t)=\Gamma(z^{-1};p,q,t).$

A solution of the key equation \eqref{ellgammaeq} in the domain $|q|>1$ is
easily found to be
$$
\gamma(z)=\frac{1}{\Gamma(q^{-1}z;p,q^{-1})}=\Gamma(pz^{-1};p,q^{-1}).
$$

As to the regime $|q|=1$, one has to abandon meromorphicity of
solutions of \eqref{ellgammaeq} in $z$ and look for an analytical
function of $u$ solving the finite-difference equation
\beq
  f(u+\omega_1) = \theta(e^{2 \pi {i} u/\omega_2};p) f(u)
\lab{eqegamma}\ee
valid for $\omega_1/\omega_2\in\mathbb{R}$.
The function $f(u)=\Gamma(e^{2\pi{i} u/\omega_2};p,q)$
solving this equation for $|q|<1$ satisfies two more equations
$$
f(u+\omega_2) = f(u), \qquad
f(u+\omega_3) = \theta(e^{2 \pi {i} u/\omega_2};q) f(u).
$$
For incommensurate $\omega_i$ these three equations define function $f(u)$ uniquely
up to multiplication by a constant. This follows from the Jacobi
theorem stating that nontrivial functions cannot have three
incommensurate periods.

For $|q|<1$ the general solution of \eqref{eqegamma} has the form $f(u)\varphi(u)$
with arbitrary periodic function $\varphi(u+\omega_1)=\varphi(u)$.
It appears that for a special choice of $\varphi(u)$ this product
defines an analytic function of $u$ even for $\omega_1/\omega_2>0$.
Such a choice has been found in \cite{spi:theta2}, where
the following modified elliptic gamma function was introduced:
\beq
G(u;{\bf \omega}):= \Gamma(e^{2 \pi {i}
u/\omega_2};p,q) \Gamma(r e^{-2 \pi {i}
u/\omega_1};\widetilde{q},r).
\lab{modeg}\ee
This function satisfies \eqref{eqegamma} and two other equations
\beq
G(u+\omega_2) = \theta(e^{2\pi i u/\omega_1};r) G(u), \quad
\lab{modeq}\ee
\beq
G(u+\omega_3) =
\frac{ \theta(e^{2 \pi {i} u/\omega_2};q)}
{ \theta(e^{-2 \pi {i} u/\omega_1};\tilde q)}G(u)
=e^{-\pi i B_{2,2}(u;\omega_1,\omega_2)} G(u),
\lab{modeq3}\ee
where $B_{2,2}$ is a second order Bernoulli polynomial
$$
B_{2,2}(u;\omega_1,\omega_2)=\frac{u^2}{\omega_1\omega_2}-\frac{u}{\omega_1}
-\frac{u}{\omega_2}+\frac{\omega_1}{6\omega_2}+\frac{\omega_2}{6\omega_1}+\frac{1}{2}.
$$
Here the exponential multiplier in \eqref{modeq3}  emerges from the
modular transformation law for the theta function
\begin{equation}
\theta\left(e^{-2\pi{i}\frac{u}{\omega_1}};
e^{-2\pi{i}\frac{\omega_2}{\omega_1}}\right)
=e^{\pi{i}B_{2,2}(u;\mathbf{\omega})} \theta\left(e^{2\pi{i}\frac{u}{\omega_2}};
e^{2\pi{i}\frac{\omega_1}{\omega_2}}\right).
\label{modthetashort}\end{equation}

{\em Exercise:} derive this relation from the modular transformation
laws for the Jacobi $\theta_1$-function
\begin{equation}
\theta_1(u/\tau|-1/\tau)
=-i\sqrt{-i\tau}\; e^{\pi iu^2/\tau}\theta_1(u|\tau)
\label{modulartheta}\end{equation}
and the Dedekind $\eta$-function
\begin{equation}
\eta(-1/\tau)=(-i\tau)^{1/2} \eta(\tau),
\quad \eta(\tau)=e^{\frac{\pi i\tau}{12}}
\left(e^{2\pi i\tau};e^{2\pi i\tau}\right)_\infty.
\label{ded}\end{equation}

Now one can check that the same three equations \eqref{eqegamma}, \eqref{modeq} and \eqref{modeq3}
and the normalization condition $G(\sum_{k=1}^3\omega_k/2;{\bf \omega})=1$ are satisfied
by the following function
\beq
G(u;\mathbf{\omega})
= e^{-\frac{\pi {i}}{3}B_{3,3}(u;\mathbf{\omega})}
\Gamma(e^{-2\pi{i} \frac{u}{\omega_3}};\tilde r,\tilde p),
\lab{modeg'}\ee
where $|\tilde p|,|\tilde r|<1$, and  $B_{3,3}(u;\mathbf{\omega})$
is the third order Bernoulli polynomial
$$
B_{3,3}(u;\mathbf{\omega})
=\frac{(u-\sum_{m=1}^3\frac{\omega_m}{2})((u-\sum_{m=1}^3\frac{\omega_m}{2})^2-\frac{1}{4}\sum_{m=1}^3\omega_m^2)}{\omega_1\omega_2\omega_3}.
$$
Since the solution of this set of equations is unique (from the nonexistence
of triply periodic functions and given normalization), we conclude that
the functions \eqref{modeg} and \eqref{modeg'} coincide.

However, from expression \eqref{modeg'}, the function
$G(u;\mathbf{\omega})$ is seen to remain a well-defined
meromorphic function of $u$ even for $\omega_1/\omega_2>0$.
Indeed, if the latter ratio is real, one can take both
 $\omega_1$ and $\omega_2$ real (since the parameters
enter only in ratios). Then one will have simultaneously
$|\tilde r|<1$ and $|\tilde p|<1$ guaranteeing convergence
of infinite products in \eqref{modeg'} only if $\omega_1/\omega_2>0$,
which gives $|q|=1$. The equality of  \eqref{modeg} and \eqref{modeg'}
is directly related to a special modular transformation
for the elliptic gamma function from the $\textrm{SL}(3,\Z)$-group
\cite{fel-var:elliptic}.

The function $G(u;{\bf \omega})$ satisfies the reflection relation
$G(a;{\bf \omega})G(b;{\bf \omega})=1,$ $a+b=\sum_{k=1}^3\omega_k$.
From \eqref{modeg'} it is not difficult to see the symmetry
$G(u;\omega_1,\omega_2,\omega_3)=G(u;\omega_2,\omega_1,\omega_3)$.

The multiple Bernoulli polynomials described above are generated
by the following expansion:
$$
\frac{x^m e^{xu}}{\prod_{k=1}^m(e^{\omega_k x}-1)}
=\sum_{n=0}^\infty B_{m,n}(u;\omega_1,\ldots,\omega_m)\frac{x^n}{n!}
$$
emerging in the theory of Barnes multiple gamma function \cite{bar:multiple}.

Let us take the limit Im$(\omega_3)\to+\infty$
and assume that Re$(\omega_1)$, Re$(\omega_2)>0$.
Then  Im$(\omega_3/\omega_1)$, Im$(\omega_3/\omega_2)\to +\infty$
and $p,r\to0$.  As a result, the expression \eqref{modeg} reduces to
$$
\stackreb{G(u;{\bf \omega})}{p,r\to 0}=\frac{(e^{2\pi iu/\omega_1}\tilde q; \tilde q)_\infty}
{(e^{2\pi i u/\omega_2}; q)_\infty}.
$$
From the representation \eqref{modeg'} one obtains a singular relation
$$
\stackreb{G(u;{\bf \omega})}{p,r\to 0}=e^{\frac{\pi i}{2}B_{2,2}(u,\omega_1,\omega_2)}
\lim_{\text{Im}(\frac{\omega_3}{\omega_1}),
\text{Im}(\frac{\omega_3}{\omega_2})\to +\infty} e^{-\pi
i\omega_3 \frac{2u-\omega_1-\omega_2}{12\omega_1\omega_2} }
\Gamma(e^{-2\pi i\frac{u}{\omega_3}};
{e^{-2\pi i\frac{\omega_1}{\omega_3}}, e^{-2\pi i\frac{\omega_2}{\omega_3}}}).
$$
For Re$(\omega_1)$, Re$(\omega_2)>0$ and $\omega_3 \to +i\infty$
this result can be rewritten as an asymptotic relation
\beq
\Gamma(e^{-2\pi v u}; e^{-2\pi v\omega_1}, e^{-2\pi v\omega_2})
\stackreb{=}{v\to 0^+}
e^{-\pi \frac{2u-\omega_1-\omega_2}{12v\omega_1\omega_2} }
\gamma^{(2)}(u;\omega_1,\omega_2),
\label{limegf}\ee
where
\beq
\gamma^{(2)}(u;\omega_1,\omega_2)
:=e^{-\frac{\pi i}{2}B_{2,2}(u;\omega_1,\omega_2)}
\frac{(e^{2\pi iu/\omega_1}\tilde q; \tilde q)_\infty}
{(e^{2\pi i u/\omega_2}; q)_\infty},
\lab{hgf}\ee
is the standard hyperbolic gamma function.

{\em Exercise:} derive the infinite product  representation \eqref{hgf}
from the integral representation
\beq
\gamma^{(2)}(u;\omega_1,\omega_2)= \exp\left(
-\textrm{p.v.}\int_\mathbb{R} \frac{e^{ux}}{(e^{\omega_1x}-1)(e^{\omega_2x}-1)}\frac{dx}{x}\right)
\lab{hgfint}\ee
with appropriate restrictions on the parameters needed for convergence.
Here ``p.v." means ``principal value", i.e. an average of integrals with the contours
passing infinitesimally above and below the singular point $x=0$.

In particular, note that for $\text{Re}(\omega_1),\text{Re}(\omega_2)>0$ and
$0<\text{Re} (u)<\text{Re}(\omega_1)+\text{Re}(\omega_2)$ integral in
\eqref{hgfint} converges and  defines $\gamma^{(2)}(u;\omega_1,\omega_2)$
as an analytic function of $u$ even for $\text{Im}(\omega_1/\omega_2)=0$,
when $|q|=1$.
The limiting relation \eqref{limegf} was rigorously established first in a different way
by Ruijsenaars \cite{rui:first}. Its uniformity was proven by Rains in \cite{rai:limits}.
The hyperbolic gamma function plays a
crucial role in the construction of $q$-hypergeometric functions in the
regime $|q|=1$ \cite{jm}. It was introduced in quantum field theory by Faddeev under the
name modular (or noncompact) quantum dilogarithm \cite{Fad94,FKV}. In a similar sense, the
elliptic gamma function has a meaning of a ``quantum"
deformation of the elliptic dilogarithm function.

\section{The elliptic beta integral}

One of the differences from ordinary hypergeometric functions and their $q$-deformations
consists in the fact that it is not straightforward to construct an equation
which is satisfied by the general elliptic hypergeometric function \eqref{EHI}.
In order to find elliptic analogues of the relations described in the introduction
one has to impose additional structural constraints on the corresponding parameters.
A basic germ, a kind of the cornerstone for building constructive
identities for such integrals is provided by
the evaluation of univariate elliptic beta integral \cite{spi:umn}.

Let complex parameters $p, q, t_j,\, j=1,\ldots, 6,$
satisfy the constraints $|p|, |q|, |t_j|<1$ and the balancing condition
$$
\prod_{j=1}^6t_j=pq.
$$
Then the following integral identity holds true
\beq
\kappa\int_\T\frac{\prod_{j=1}^6
\eg(t_jx^{{\pm 1}};p,q)}{\eg(x^{\pm 2};p,q)}\frac{dx}{x}
=\prod_{1\leq j<k\leq6}\eg(t_jt_k;p,q),
\lab{elbeta}\ee
where $\T$ is the unit circle of positive orientation and
$$
\kappa=\frac{(p;p)_\infty(q;q)_\infty}{4\pi {i}}.
$$

We sketch the proof of this statement suggested in \cite{spi:short}.
Note first that the integrand has poles at the points
$z=t_jq^ap^b,\, j=1,\ldots,6, a,b\in\Z_{\geq0}$, converging to zero,
and their reciprocals $z=t_j^{-1}q^{-a}p^{-b}$, diverging to infinity.
The integration contour $\T$ separates these sets of poles.

Now we apply the gamma function inversion
$$
\Gamma(t_6x;p,q)=\frac{1}{\Gamma(pq/(t_6 x);p,q)}=\frac{1}{\Gamma(Ax^{-1};p,q)},
\quad  A:=\prod_{m=1}^5t_m,
$$
and rewrite the integral evaluation as
$$
I(t_1,\ldots, t_5;p,q)=\kappa\int_{\T}\rho(x;t_1,\dots,t_5;p,q)\frac{dx}{x}=1,
$$
where
$$
\rho(x;t_1,\dots,t_5;p,q)
=\frac{\prod_{j=1}^5
\Gamma(t_jx^{\pm1},t_j^{-1}A;p,q)}
{\Gamma(x^{\pm2},Ax^{\pm1};p,q)
\prod_{1\le i < j\le 5}\Gamma(t_it_j;p,q)}.
$$

This kernel function satisfies the $q$-difference equation
\beq
\rho(x;qt_1)-\rho(x;t_1)
=g(q^{-1}x)\rho(q^{-1}x;t_1)-g(x)\rho(x;t_1)
\lab{keyeq}\ee
with
$$
g(x)=
\frac{\prod_{m=1}^5\theta(t_mx;p)}{\prod_{m=2}^5\theta(t_1t_m;p)}
\frac{\theta(t_1A;p)}
{\theta(x^2,xA;p)}\frac{t_1}{x}.
$$
Dividing \eqref{keyeq} by the $\rho$-function, one comes to
the following elliptic functions identity
\begin{eqnarray*} && \makebox[-1em]{}
\frac{\theta(t_1x,t_1x^{-1};p)}{\theta(Ax,Ax^{-1};p)}
\prod_{m=2}^5\frac{\theta(At_m^{-1};p)}{\theta(t_1t_m;p)}-1 =
\frac{t_1\theta(t_1A;p)}{x\theta(x^2;p)\prod_{m=2}^5\theta(t_1t_m;p)}
\\ &&  \makebox[2em]{}
\times
\left(\frac{x^4\prod_{m=1}^5\theta(t_mx^{-1};p)}{\theta(Ax^{-1};p)} -
\frac{\prod_{m=1}^5\theta(t_mx;p)}{\theta(Ax;p)}\right).
\end{eqnarray*}
which we suggest to prove as an exercise (compare the poles in $x$ and their residues in the
parallelogram of periods of the left- and right-hand side expressions and verify the
identity for a particular value of $x$).

Integrating equation \eqref{keyeq} over $x\in\T$ one obtains the
relation
$$
I(qt_1)-I(t_1)=\Big(\int_{q^{-1}\T} -\int_\T\Big)g(x)\rho(x;t_1)\frac{dx}{x}.
$$
Consideration of the poles of the function $g(x)\rho(x;t_1)$ shows
that for $|t_1|< |q|$ it does not have singularities inside the annulus
bounded by $\T$ and the circle of radius $|q|^{-1}$ denoted as $q^{-1}\T$.
As a result, the right-hand side expression in the above relation
vanishes and the equality $I(qt_1)=I(t_1)$ emerges in a natural way.
After permuting $p$ and $q$
in the above considerations and imposing the additional
constraint $|t_1|<|p|$, it becomes possible to write $I(pt_1)=I(t_1)$.
Now, the Jacobi theorem on the absence of periodic functions with three
incommensurate periods and $t_k$-permutational symmetry show
that $I$ does not depend on parameters, $I=I(p,q)$.

In order to compute this constant, one can consider the limit $t_1t_2\to 1$
when two pairs of residues pinch the contour of integration.
After crossing a pair of poles and picking up the residues,
one can see that the integral part vanishes in this limit, and the
contribution of residues sums exactly to the needed value $I=1$.

The derived elliptic beta integral evaluation represents a unique
relation due to the following facts.
First of all, it represents an elliptic extension of Newton's binomial theorem
$_1F_0(a;x)=(1-x)^{-a}$ and its $q$-analogue
$$
_1\varphi_0(t;q;x)=\sum_{n=0}^\infty \frac{(t;q)_n}{(q;q)_n}x^n
=\frac{(tx;q)_\infty}{(x;q)_\infty}, \quad
(t;q)_n=\prod_{k=0}^{n-1}(1-tq^k).
$$
After setting $t=q^a$, for $q\to 1^-$ one has $_1\varphi_0(q^a;q;x)\to {}_1F_0(a;x)$.
This yields the useful relation
\beq
\lim_{q\to 1^-} \frac{(q^ax;q)_\infty}{(x;q)_\infty} =(1-x)^{-a}.
\lab{limbin}\ee

As to the terminating series version of the binomial theorem, its elliptic analogue is
given by the Frenkel-Turaev sum \cite{FT}, which can be obtained by a reduction of  \eqref{elbeta}.
To derive this sum let us take the limit $t_4t_5\to q^{-N}$ for some positive integer $N$.
More precisely, let us take parameter $t_5$ from inside $\T$ to outside such
that $|pt_5|, |q^{N+1}t_5| <1< |t_5|$ and keep other parameters inside $\T$ in
generic positions.
Formula \eqref{elbeta} will remain intact if we replace the contour $\T$ by $C$ which separates
sequences of poles converging to zero from their reciprocals. Consideration of the poles
related to the parameters $t_4$ and $t_5$ shows that if $t_4t_5\to q^{-N}$ then
$4(N+1)$ poles start to pinch pairwise two parts of the contour $C$ lying outside and inside $\T$.
As a result both, the left- and right-hand side expressions in \eqref{elbeta} start to diverge.

To compute the limiting formula, resolve the balancing condition $t_6=pq/\prod_{k=1}^5t_k$
and denote
$\rho_E(z,\underline{t})=\prod_{m=1}^5\eg(t_mz^{\pm 1};p,q)/\eg(z^{\pm2},\prod_{k=1}^5t_k z^{\pm 1};p,q)$.
Let us force the contour $C$ to cross $2(N+1)$ poles $z=(t_5q^k)^{\pm1}$, $k=0,\ldots,N$.
Then the Cauchy theorem states that:
\begin{eqnarray}  \nonumber &&
\kappa\int_C \rho_E(z,\underline{t})\frac{d z}{z} =
\kappa\int_\mathbb{T} \rho_E(z,\underline{t})\frac{d z}{z}
\\  && \makebox[2em]{}
+\frac{\prod_{m=1}^4\eg(t_mt_5^{\pm 1};p,q)}
{\eg(t_5^{-2},\prod_{k=1}^5t_kt_5^{\pm 1};p,q)}\sum_{n=0}^N
\frac{\theta(t_5^2q^{2n};p)}{\theta(t_5^2;p)}
\prod_{m=0}^5 \frac{\theta(t_mt_5)_n} {\theta(qt_m^{-1}t_5)_n}\, q^n.
\lab{res}\end{eqnarray}
Here we denoted $t_0=q/\prod_{m=1}^5t_m=t_6/p$ and used the elliptic Pochhammer symbol
$$
\theta(t)_n=\prod_{j=0}^{n-1}\theta(tq^j;p)=\frac{\eg(tq^n;p,q)}{\eg(t;p,q)}, \qquad
\theta(t_1,\ldots,t_k)_n:=\prod_{j=1}^k\theta(t_j)_n.
$$
The residues are computed using the limiting relation \eqref{limit}.

Now we take the desired limit $t_5t_4\to q^{-N}$.  The integral over the unit circle $\mathbb{T}$
stays finite, since the integrand is nonsingular on $\T$, whereas the sum of residues and the
value of the original integral diverge. Dividing expression \eqref{res} and its evaluation
 \eqref{elbeta} by $\Gamma(t_4t_5;p,q)$, for $t_5t_4= q^{-N}$ one obtains the Frenkel-Turaev sum
\begin{equation}
{}_{10}V_9(t_5^2;t_0t_5,t_1t_5,t_2t_5,t_3t_5,q^{-N})=
\frac{\theta(qt_5^2,\frac{q}{t_1t_2},
\frac{q}{t_1t_3},\frac{q}{t_2t_3})_N }
{\theta(\frac{q}{t_1t_2t_3t_5},\frac{qt_5}{t_1},
\frac{qt_5}{t_2},\frac{qt_5}{t_3})_N}.
\label{ft-sum}\end{equation}
Here we use general notation for the very-well poised elliptic hypergeometric series
introduced in \cite{spi:theta1}
\begin{equation}
{}_{m+1}V_{m}(t_0;t_1,\ldots,t_{m-4};q,p)
=\sum_{n=0}^\infty \frac{\theta(t_0q^{2n};p)}{\theta(t_0;p)}\prod_{k=0}^{m-4}
\frac{\theta(t_k)_n}{\theta(qt_0t_k^{-1})_n}q^n
\label{Vseries}\end{equation}
with the balancing condition $\prod_{k=1}^{m-4}t_k=t_0^{(m-5)/2}q^{(m-7)/2}$
and the assumption that the series terminates because one of the parameters has the form
$t_j=q^{-N}$. For $p\to 0$ the series \eqref{Vseries} with fixed parameters
reduces to the very-well poised balanced
${}_{m-1}\varphi_{m-2}$ series \cite{gas-rah:basic}.
The original derivation of \eqref{ft-sum} in \cite{FT} is completely different from the given one
which was suggested in \cite{die-spi:elliptic}.  Multivariable extensions of the
elliptic hypergeometric series were considered for the first time in \cite{Warnaar2002}.

{\em Exercise:} verify the above derivation of \eqref{ft-sum} by completing all the details.

The next important property of the integral \eqref{elbeta} is that it represents
the top known generalization of the Euler beta integral
\eqref{beta}. In particular, in the limit $p \to 0$, taken for fixed $t_1,\ldots,t_5$ and
$t_6 \propto p$, one obtains the Rahman $q$-beta integral \cite{rah:integral}.
Subsequent turning one of the parameters to zero yields the Askey-Wilson $q$-beta
integral whose reduction to \eqref{beta} was explicitly described in \cite{gas-rah:basic}.
More complicated degenerations of the elliptic beta integral are considered in \cite{SarSpi2}.

Integral \eqref{elbeta} serves as the measure for a biorthogonality relation of
specific two-index functions, defined as products of two $_{12}V_{11}$ series,
which generalize the Askey-Wilson, Jacobi and other
classical orthogonal polynomials \cite{spi:theta2}. These functions comprise also
Rahman's continuous biorthogonal rational functions \cite{rah:integral}.
The discrete measure analogues of these functions based on the Frenkel-Turaev sum
were defined in \cite{spi-zhe:spectral}.

Integral \eqref{elbeta} is a germ for constructing
infinitely many elliptic hypergeometric integrals admitting
exact evaluation. It generates an elliptic Fourier transform \cite{spi:bailey2,spi-war:inversions},
associated with an integral generalization of the Bailey chains techniques \cite{Warnaar},
integral operator realization of Coxeter relations \cite{DS}, the star-triangle relation \cite{BS}, and
the Yang-Baxter equation. Identity \eqref{elbeta} emerged in four dimensional
supersymmetric quantum field theory as an equality of superconformal indices
of two specific models \cite{DO}. Some of these unique features of the elliptic
beta integral are described in more detail in the following.

The very first proof of formula \eqref{elbeta} was obtained
using the contiguous relation  for integrals \eqref{cont-1} and expansion in small $p$
\cite{spi:umn}, when the limiting $t_j=0$ points enter the
domain of analyticity of the expansion coefficients. A further
refinement of such expansion arguments was suggested in \cite{rai:sigma},
when the equality of formal series in $p$ in the left- and right-hand sides
of identities is reached by establishing their
rationality and coincidence on an infinite discrete set of parameter values.
This gives another proof of the above formula based on the theory of Askey-Wilson
polynomials and Frenkel-Turaev sum. The proof of \cite{spi:short} given above is self-contained -
it does not require knowledge of any system of orthogonal functions and
uses only a simple elliptic function identity.

\section{An elliptic extension of the Euler-Gauss hypergeometric function}

There are many generalizations of the $F(a,b;c;x)$ hypergeometric function.
Let us describe the one related to the elliptic beta integral in a
way as the beta function \eqref{beta} is connected to \eqref{2F1}.
It is necessary to take two base variables $p, q$, $|p|,|q|<1$,
and eight parameters $t_1,\ldots, t_8\in\C^\times  $ satisfying the
balancing condition
$
\prod_{j=1}^8 t_j=p^2q^2.
$
Then, under additional constraints $|t_j|<1$, an elliptic analogue of
the Euler-Gauss hy\-per\-geo\-met\-ric function is defined by
the integral \cite{spi:theta2}
\be
V(\underline{t})\equiv V(t_1,\ldots,t_8;p,q):=\kappa\int_\T\frac{\prod_{j=1}^8
\eg(t_jz^{\pm 1};p,q)}{\eg(z^{\pm 2};p,q)}\frac{dz}{z}.
\lab{ehf}\ee
By deforming the integration contour it is possible to partially
relax the constraints on the parameters. Analytic continuation of
\eqref{ehf} is achieved by increasing the absolute values of
parameters and computing the residues of the integrand
poles, so that the analytically continued function becomes a sum
of the integral over some fixed contour and residues of the
poles crossed by the contour. From this procedure one
can see that the $V$-function is meromorphic
for all values of parameters $t_j\in\C^\times $, when the contour of integration
is not pinched which may happen for $t_jt_k=q^{-a}p^{-b},\, a,b\in\Z_{\geq0}$.
It appears that the potential singularities from $t_j^2=q^{-a}p^{-b}$ do not
contribute and the product $\prod_{1\leq j<k\leq 8}(t_jt_k;p,q)_\infty V(\underline{t})$
becomes a holomorphic function of the parameters \cite{rai:trans}.
For  particular values of the  parameters $t_j$ the $V$-function has
delta-function type singularities  \cite{S5,SV3}.
We remark that the expression \eqref{ehf} can be reduced to both Euler and Barnes type integral
representations for $F(a,b;c;x)$.

The hypergeometric function  \eqref{eulerint} is reduced to Euler's beta integral
for $a=0$. In a similar way, its elliptic counterpart \eqref{ehf} reduces to
the elliptic beta integral if a pair of parameters is constrained as $t_jt_k=pq$, $j\neq k$,
as follows from the inversion relation \eqref{inv}.

Consider now symmetry transformations for the $V$-function.
An evident symmetry is the possibility to permute bases $p$ and $q$.
For describing symmetries in the parameters we remind some simplest facts
from the theory of root systems and corresponding Weyl groups.
Consider $\R^n$ with an orthonormal basis $e_i\in \R^n$,
 $\langle e_i, e_j\rangle=\delta_{ij}$.
For any $x\in\R^n$ define its reflection with respect to the hyperplane orthogonal to
some $v\in\R^n$:
$$
x\to S_v(x)=x-\frac{ 2\langle v, x\rangle}{ \langle v, v\rangle}\,v, \qquad S^2_v=1.
$$
If $x=const\cdot v$, then $S_v(x)= -x$. For $\langle v, x\rangle=0$ one has $S_v(x)=x$.

Define $R$ as some set of vectors $\alpha_1,\ldots, \alpha_m\in \R^n$,
forming a basis. If for any $\alpha, \beta\in R$, $S_\alpha(\beta)\in R$,
then $R$ is called a root system. The reflections $W=\{S_\alpha\}$ form
a finite subgroup of the rotation group $O(n)$.
The vectors $\alpha_j$ are called the roots and the dimensionality of the
space where they are defined is the rank of the root system.

If for all $\alpha, \beta\in R$ one has the integrality
$
2\langle \alpha, \beta\rangle/\langle \alpha, \alpha\rangle\in\Z,
$
then $R$ is called the crystallographic root system and $W$ -- the Weyl group.

If the only multiples of a root $\alpha$ in  $R$  are $\pm\alpha$
then $R$ is called reduced and it is known to be related to a
semi-simple Lie algebra. For such cases there exist four
irreducible (i.e. indecomposable to direct sums)
infinite classical series of root systems:
$A_n, B_n, C_n$, $D_n$, and five exceptional cases:
$G_2, F_4, E_6, E_7, E_8$.

Let us describe a few examples of root systems used in the following.

1) $A_n$ system: take $\mathbb{E} \in \R^{n+1}$ orthogonal to
$\sum_{j=1}^{n+1}e_j$, i.e. for $u=\sum_{j=1}^{n+1}u_je_j \in \mathbb{E}$
one has $\sum_{j=1}^{n+1}u_j=0$.
Then $R_{A_n}=\{e_i-e_j,\, i\neq j\}\big|_{i,j=1,\ldots,n+1}$
and the Weyl group is the permutation group $W(A_n)=S_{n+1}$.

2) $C_n$ system: take in $\R^n$ the roots
$R_{C_n}=\{\pm 2e_i; \pm e_i \pm e_j,\, i< j\}\big|_{i,j=1,\ldots,n}$
and $W=S_{n}\times \Z_2^n$. The only non-reduced root system is $R_{BC_n}$,
which contains the roots of the $C_n$ system and additionally
$\{\pm e_1,\ldots, \pm e_n\}$.

3) $E_7$ system: take $\mathbb{E} \in \R^{8}$ orthogonal to
$\sum_{j=1}^{8}e_j$, as for the $A_7$ root system.
Then $R_{E_7}=\{e_j-e_k,\, j\neq k; \frac{1}{2}\sum_{l=1}^8\mu_j e_j,\,
\mu_j=\pm 1 \text{ with four values } \mu_j=1 \}
\big|_{j,k=1,\ldots,8}$, and $W(E_7)$ is a particular finite group of order $72\cdot 8!$

As we will see, elliptic hypergeometric integrals are naturally related to
the root systems in two qualitatively different ways.

The $V$-function is evidently invariant under the $S_8$-group of
permutations of parameters $t_j$. It is the Weyl group of the $A_7$ root
system. Consider now the double integral
$$
\kappa \int_{\T^2}
\frac{\prod_{j=1}^4\eg(a_jz^{\pm 1}, b_jw^{\pm 1};p,q)\;\eg(cz^{\pm 1} w^{\pm 1};p,q)}
{\eg(z^{\pm 2},w^{\pm 2};p,q)} \frac{dz}{z}\frac{dw}{w},
$$
where complex parameters $a_j,b_j,c\in \C^\times  $ are constrained as
$|a_j|,|b_j|,|c|<1$ and satisfy the balancing conditions
$$
c^2\prod_{j=1}^4a_j= c^2\prod_{j=1}^4b_j=pq.
$$
Since we integrate over compact domains, the order of integrations does
not matter. The integrals over $z$ or $w$ are separately computable
due to the key formula \eqref{elbeta}. Taking these integrals in the
different order we come to the following transformation formula:
\begin{eqnarray*}
&& V(ca_1,\ldots, ca_4,b_1,\ldots,b_4)
=\prod_{1\le j<k\le 4}\frac{\Gamma(b_jb_k;p,q)}{\Gamma(a_ja_k;p,q)}V(a_1,\ldots, a_4,cb_1,\ldots,cb_4),
\end{eqnarray*}
which can be rewritten in a more symmetric form
\beq
V(\underline{t})=\prod_{1\le j<k\le 4}\eg(t_jt_k,t_{j+4}t_{k+4};p,q)\,
V(\underline{s}),
\lab{E7-1}\ee
where the parameters $t_j$ and $s_j$ are related to each other as
$$
\left\{
\begin{array}{cl}
s_j =\sqrt{\frac{pq}{t_1t_2t_3t_4}}\, t_j,&   j=1,2,3,4  \\
s_j = \sqrt{\frac{pq}{t_5t_6t_7t_8} }\, t_j, &    j=5,6,7,8
\end{array}
\right.
$$
and satisfy the constraints $|t_j|, |s_j|<1$
matching the integration contour $\T$ on both sides of \eqref{E7-1}.

The function $V(\underline{t})$ appeared for the first time during
the derivation of this fundamental relation in  \cite{spi:theta2}.
Let us write $t_j=e^{2\pi{i} x_j}(pq)^{1/4}$
and $s_j=e^{2\pi{i} y_j}(pq)^{1/4}$, $j=1,\ldots,8$.
From the balancing condition
we find $\sum_{j=1}^8 x_j=\sum_{j=1}^8 y_j=0$. Now it is not
difficult to see that the transformation of parameters in \eqref{E7-1} is equivalent
to the relation $y_j=x_j-\frac{\mu}{4}\sum_{k=1}^4(x_k-x_{k+4})$
with $\mu=1$ for $j=1,\ldots,4$ and $\mu=-1$ for $j=5,\ldots,8$,
which precisely corresponds to the reflection $y=S_v(x)$ for the vector
$v=(\sum_{k=1}^4e_i-\sum_{k=5}^8 e_i)/2$ of the canonical length
$\langle v, v\rangle=2$. Permuting in $\left(8 \atop 4\right)=70$ nontrivial ways
the basis vectors in this $v$ one comes to the roots of the exceptional root
system $E_7$ extending the $A_7$ root system, as described above.

Now one can consider all admissible $W(E_7)$ Weyl group reflections
acting on the $V$-function. For instance, it is possible to repeat
reflection \eqref{E7-1} for the second time using the
root $v=(e_3+e_4+e_5+e_6-e_1-e_2-e_7-e_8)/2$ and permute in the resulting
relation parameters in all possible ways. This yields the following
symmetry transformation
\beq
V(\underline{t})=\prod_{j,k=1}^4
\eg(t_jt_{k+4};p,q)\ V(T^{1\over 2}\!/t_1,\ldots,T^{1\over 2}\!/t_4,
S^{1\over 2}\!/t_5,\ldots, S^{1\over 2}\!/t_8),
\lab{E7-2}\ee
where $ T=t_1t_2t_3t_4$, $ S=t_5t_6t_7t_8$ and one has the constraints
$|T|^{1/2}<|t_j|<1,$ $|S|^{1/2}<|t_{j+4}|<1,\, j=1,2,3,4$, in order
to have $\T$ as the integration contour on both sides.

Finally, let us equate expressions on the right-hand sides of relations
\eqref{E7-1} and \eqref{E7-2}. After rewriting the resulting equality
in terms of the parameters $s_j$, it takes the form
\beq
V(\underline{s})=\prod_{1\le j<k\le 8}\eg(s_js_k;p,q)\,
V(\sqrt{pq}/s_1,\ldots,\sqrt{pq}/s_8),
\lab{E7-3}\ee
where $|pq|^{1/2}<|s_j|<1$ for all $j$. The key generating relation
\eqref{E7-1} was discovered in \cite{spi:theta2}.
Transformations \eqref{E7-2} and \eqref{E7-3} were
proved in a different way in \cite{rai:trans}, where the
identification of these transformations with the group $W(E_7)$ was made.

Although these three identities for the $V$-function have
different form, they are tied by the symmetry group. As it will be shown
later on, the multiple elliptic hypergeometric integrals have
transformations which can be considered as their separate generalizations,
i.e. different elements of $W(E_7)$ may have individual
multivariable extensions.

Let us identify parameters in \eqref{addition} as $y=t_1, x=t_2, w=t_3$ and
multiply this addition formula by
$\rho(z;\underline{t})=\prod_{j=1}^8\Gamma(t_jz^{\pm 1};p,q)/\Gamma(z^{\pm 2};p,q)$
 with the balancing condition $\prod_{j=1}^8t_j=p^2q$. Then we can write
$$
t_3\theta(t_2t_3^{\pm1};p)\rho(z;qt_1,t_2,\ldots)
+t_1\theta(t_3t_1^{\pm1};p)\rho(z;t_1,qt_2,\ldots)
+t_2\theta(t_1t_2^{\pm1};p)\rho(z;t_1,t_2,qt_3,\ldots)=0.
$$
Integrating this relation over $z\in\T$ we obtain
 the following contiguous relation
\begin{equation}
\frac{t_1V(qt_1)}{\theta(t_1t_2^{\pm1},t_1t_3^{\pm 1};p)}
+\frac{t_2V(qt_2)}{\theta(t_2t_1^{\pm1},t_2t_3^{\pm 1};p)}
+\frac{t_3V(qt_3)}{\theta(t_3t_1^{\pm1},t_3t_2^{\pm 1};p)} = 0,
\label{cont-1}\end{equation}
where $V(qt_j)$ denotes the $V(\underline{t};p,q)$-function
with the parameter $t_j$ replaced by $qt_j$, so that the
balancing condition takes the form indicated above.

Applying symmetry relations discussed in the previous section to
the $V$-functions in \eqref{cont-1} one obtains many differently looking identities.
In particular, substitution of the third transformation \eqref{E7-3}
yields the contiguous relation
\begin{equation}
\frac{\prod_{j=4}^8\theta\left(t_1t_j/q;p\right)V(t_1/q)}
{t_1\theta(t_2/t_1,t_3/t_1;p)}
+\frac{\prod_{j=4}^8\theta\left(t_2t_j/q;p\right)V(t_2/q)}
{t_2\theta(t_1/t_2,t_3/t_2;p)}
+\frac{\prod_{j=4}^8\theta\left(t_3t_j/q;p\right)V(t_3/q)}
{t_3\theta(t_1/t_3,t_2/t_3;p)} = 0,
\label{cont-3}\end{equation}
where $\prod_{j=1}^8t_j=p^2q^3$.

 Consider now three equations:
1) the equation obtained from \eqref{cont-1} after the replacement
$t_1 \to q^{-1}t_1$, 2) the one obtained from \eqref{cont-1}
after the replacement $t_2 \to q^{-1}t_2$, and 3) the $t_3\to qt_3$
transformed version of \eqref{cont-3}.
Eliminating from them the functions $V(q^{-1}t_1,qt_3)$ and $V(q^{-1}t_2,qt_3)$
we come to the elliptic hypergeometric equation \cite{spi:cs}:
\begin{eqnarray}\label{eheq1}
&& \makebox[-2em]{}
\mathcal{A}(t_1,t_2,\ldots,t_8,q;p)\Big(U(qt_1,q^{-1}t_2;p,q)-U(\underline{t};p,q)\Big)
\\ &&
+\mathcal{A}(t_2,t_1,\ldots,t_8,q;p)\Big(U(q^{-1}t_1,qt_2,;p,q)-U(\underline{t};p,q)\Big)
+ U(\underline{t};p,q)=0,
\nonumber\end{eqnarray}
where
\begin{equation}
 \mathcal{A}(t_1,\ldots, t_8,q;p)=\frac{\theta(t_1/qt_3,t_3t_1,t_3/t_1;p)}
                 {\theta(t_1/t_2,t_2/qt_1,t_1t_2/q;p)}
\prod_{k=4}^8\frac{\theta(t_2t_k/q;p)}{\theta(t_3t_k;p)}
\label{Apot}\end{equation}
and
$$
U(\underline{t}; p,q):=\frac{V(\underline{t};p,q)}
{\prod_{k=1}^2 \eg(t_kt_3^{\pm 1};p,q)}.
$$

{\em Exercise:} verify that the coefficient $\mathcal{A}(t_1,\ldots, t_8,q;p)$ is invariant under the
transformations $t_j\to p^{n_j}t_j,\, q\to p^nq$ for any $n_j, n\in\Z$
preserving the balancing condition $\prod_{j=1}^8t_j=p^2q^2$.

Expressing $t_1$ in terms of $t_2$ (or vice versa) via the balancing condition,
one sees that \eqref{eheq1} is actually a second order $q$-difference equation
in $t_2$ (or $t_1$). It shows that the elliptic hypergeometric integrals may
emerge as solutions of particular finite-difference equations with elliptic
function coefficients.

Since $\mathcal{A}(p^{-1}t_1,pt_2,\ldots)=
\mathcal{A}(t_1,t_2,\ldots)$, the function $U(p^{-1}t_1,pt_2)$ defines
the second independent solution of \eqref{eheq1}.
Let us multiply \eqref{eheq1} by $U(p^{-1}t_1,pt_2)$ and the equation
for $U(p^{-1}t_1,pt_2)$ by $U(t_1,t_2)$ and subtract one from
another. This yields
\begin{equation}\label{cas-eqn}
\mathcal{A}(t_1,t_2,\ldots t_8,q;p){D}(p^{-1}t_1,q^{-1}t_2)
=\mathcal{A}(t_2,t_1,t_3,\ldots,q;p){D}(p^{-1}q^{-1}t_1,t_2),
\end{equation}
where
$$
{D}(t_1,t_2)=U(qpt_1,t_2)U(t_1,pqt_2)-U(qt_1,pt_2)U(pt_1,qt_2)
$$
is the $t_1\to pt_1$ and $t_2\to qt_2$ renormalized
version of the Casoratian (discrete Wronskian) with the balancing condition for
$U$-function parameters $\prod_{j=1}^8t_j=pq$.

Let $t_2$ be an independent variable. Then $t_1\propto 1/t_2$ due to the balancing condition.
Therefore, after denoting $f(t_2):={D}(t_1,t_2)$, relation \eqref{cas-eqn}  is nothing else than the following
first order $q$-difference equation in $t_2$:
\begin{align}\nonumber
f(qt_2)& =\frac{\mathcal{A}(pt_1,qt_2,t_3,\ldots,q;p)}
{\mathcal{A}(qt_2,pt_1,t_3,\ldots,q;p)}f(t_2)
\\
& =-\frac{t_1}{qt_2}\frac{\theta(t_1/q^2t_2,t_1/qt_3,t_1^{-1}t_3^{\pm 1};p)}
            {\theta(t_2/t_1,t_2/t_3, q^{-1}t_2^{-1}t_3^{\pm 1};p)}
\prod_{k=4}^8\frac{\theta(t_2t_k;p)}{\theta(t_1t_k/q;p)}f(t_2).
\nonumber\end{align}
Its general solution has the form
\beq
{D}(t_1,t_2)=\varphi(t_2)\;\frac{\prod_{k=3}^8\eg(t_1t_k,t_2t_k) }
{\eg(t_1/t_2,t_2/t_1)}\prod_{k=1}^2\frac{\eg(t_k^{-1}t_3^{\pm 1};p,q)}
{\eg(t_kt_3^{\pm 1};p,q)},
\lab{genD}\ee
where $\varphi(qt_2)=\varphi(t_2)$. Since ${D}(t_1,t_2)$ is symmetric in $p$ and $q$,
we can repeat the above consideration with permuted $p$ and $q$,
which yields $\varphi(pt_2)=\varphi(t_2)$. By the Jacobi theorem, for incommensurate
$p$ and $q$ this proves that $\varphi$ does not depend on $t_2$.

{\em Exercise:} compute the constant  $\varphi$ by taking the limit $t_2\to 1/t_3$ and
using the residue calculus. Show that
$$
\varphi=\frac{\prod_{3\leq j<k\leq 8}\eg(t_jt_k;p,q)}{\eg(t_1^{-1}t_2^{-1};p,q)},
$$
which yields the following quadratic relation
for the elliptic hypergeometric function \cite{RS}
\begin{equation}
V(pqt_1,t_2)V(t_1,pqt_2)-t_1^{-2}t_2^{-2} V(qt_1,pt_2)V(pt_1,qt_2)
=\frac{\prod_{1\leq j<k\leq8}\eg(t_jt_k;p,q)}{\eg(t_1^{\pm1}t_2^{\pm1};p,q)}.
\label{V-det}\end{equation}

Solutions of the elliptic hypergeometric equation  \eqref{eheq1}
which we discussed so far are defined for $|q|<1$.
However, the equation itself does not assume such a constraint.
In order to build its solutions in other domains of values of $q$
one can use symmetries of the equation \eqref{eheq1} which are not symmetries
of the described solutions. In particular, the following relation
holds true
$$
\mathcal{A}\left(\frac{p^{1/2}}{t_1}, \ldots,\frac{p^{1/2}}{t_8},q;p\right)
=\mathcal{A}\left(t_1,\ldots,t_8,q^{-1};p\right).
$$
This means that the scalings $t_j\to p^{1/2}/t_j$, $j=1,\ldots, 8,$
transform \eqref{eheq1} to the same equation with the replacement
of the base $q\to q^{-1}$. The inversion $q\to 1/q$ takes place also
if one replaces $t_j\to p^{n_j}/t_j$ with integer $n_j$,
$\sum_{j=1}^8n_j=4$. So, in the regime  $|q|>1$ one obtains the
following particular solution of \eqref{eheq1} \cite{RS}
\begin{eqnarray}
U_{|q|>1}(\underline{t};q,p)=\frac{V(p^{1/2}/t_1,\ldots,p^{1/2}/t_8;p, q^{-1})}
{\prod_{k=1}^2\Gamma(p/t_kt_3,t_3/t_k;p,q^{-1}) }.
\label{q>1}\end{eqnarray}

In order to obtain solutions of the elliptic hypergeometric equation
on the unit circle $|q|=1$, it is necessary to use the modified
elliptic gamma function $G(u;\mathbf{\omega})$. Indeed, we can replace in the definitions
of the elliptic beta integral and the $V$-function the function
$\Gamma(z;p,q)$ by $G(u;\mathbf{\omega})$ and repeat all the considerations
anew. Because the functional equations for these elliptic gamma functions
are similar, one will obtain formulas analogous to those presented
above. But from the representation \eqref{modeg'} it follows that
the difference between them lies only in the exponential factors
containing the Bernoulli polynomials. As shown in \cite{die-spi:unit},
these factors can be removed reducing everything to a modular
transformed version of the described above relations. By construction,
such relations remain well defined even if $|q|=1$. At the level
of equation \eqref{eheq1} one has the modular invariance
$$
\mathcal{A}(e^{2\pi i \frac{g_1}{\omega_2}},e^{2\pi i \frac{g_2}{\omega_2}},
\ldots e^{2\pi i \frac{g_8}{\omega_2}},e^{2\pi i \frac{\omega_1}{\omega_2}};
e^{2\pi i \frac{\omega_3}{\omega_2}})
=\mathcal{A}(\ldots)
\big|_{(\omega_2,\omega_3)\to (-\omega_3,\omega_2)}.
$$
Therefore, a solution of \eqref{eheq1} valid for $|q|=1$
is obtained by using the described parametrization of variables and
by making a particular modular transformation
$$
U_{|q|=1}(\underline{t};p,q)=
U(\underline{t};p,q)\big|_{(\omega_2,\omega_3)\to (-\omega_3,\omega_2)}.
$$

Let us give another form of the elliptic hypergeometric equation.
We single out the variable $x$ by setting $t_1=cx,\ t_2=c/x$ and denote
$$
\ve_1=\frac{c}{t_3}, \quad \ve_2=\frac{\ve_1}{q},\quad
\ve_3=ct_3p^4, \quad
\ve_k=\frac{q}{ct_k},\; k=4,\ldots,8.
$$
Since $c=\sqrt{t_1t_2}$, one has the same balancing condition
$\prod_{k=1}^8\ve_k=p^2q^2$. Evidently, scalings of parameters
of the $U$-function in \eqref{eheq1} are equivalent to the shifts
$x\to q^{\pm1} x$. After the replacement of
$U(\underline{t})$ by some unknown function $f(x)$, \eqref{eheq1}
becomes a $q$-difference equation of the second order
of the following symmetric form
\begin{eqnarray}
&& A(x)\left( f(qx)-f(x)\right)
+ A(x^{-1})\left( f(q^{-1}x)-f(x)\right) + \nu f(x)=0,
\label{eheq2}
\\ && \qquad
A(x)=\frac{\prod_{k=1}^8 \theta(\ve_kx;p)}{\theta(x^2,qx^2;p)},
\qquad
\nu=\prod_{k=3}^8\theta\left(\frac{\ve_k\ve_1}{q};p\right).
\lab{pot}\end{eqnarray}
Note that here $\ve_k$-variables are constrained not only by
the balancing condition, but also by the additional relation $\ve_2=\ve_1/q$.

Clearly equation \eqref{eheq2} has only $S_6$-symmetry in parameters $\ve_k$,
$k=3,\ldots,8$. However, as noticed by Zagier, the potential $\mathcal{A}$ from
\eqref{Apot} itself can be written  in a completely $S_8$-symmetric form.
Indeed, denote
$$
u_1=\frac{t_1}{t_3},\; u_2=\frac{t_1}{qt_3},\; u_3=\frac{1}{t_1t_3},\;
u_k=\frac{t_kt_2}{q},\; k=4,\ldots,8,\;
\lambda=\frac{t_2}{qt_3}.
$$
Then one can write
$$
\mathcal{A}(t_1,\ldots, t_8,q;p)=\frac{\lambda^2}{p^2}
\prod_{k=1}^8\frac{\theta(u_k;p)}{\theta(v_k;p)},\quad
u_kv_k=\lambda, \quad \prod_{k=1}^8 u_k=p^2\lambda^4.
$$
All $u_k$ variables are independent and $\lambda$ is determined by
their product, i.e. the presence of the $S_8$ symmetry becomes evident.

Because of the distinguished role of the elliptic hypergeometric equation
it is interesting to know all its roots of origin.
It appears \cite{spi:cs} that equation \eqref{eheq2} is related to the eigenvalue problem $H\psi=E\psi$
for the restricted one particle Hamiltonian of the van Diejen model \cite{die:integrability}.
Namely, one has to take special eigenvalue $E=-\nu$
and impose two additional constraints on the parameters of the general model ---
the balancing condition and $\ve_2=\ve_1/q$.
Another place where this equation emerges in a natural way is the
theory of elliptic Painlev\'e equation \cite{sak}. Namely, for a special
restriction on the geometry of this equation it linearizes exactly to the
elliptic hypergeometric equation \cite{kmnoy}. In a related subject
it emerges as the simplest rigid equation in the elliptic isomonodromy problem
\cite{rai:painleve,rai:noncom}. A list of degenerations of the $V$-function
to the lower level hypergeometric functions is considered in detail in \cite{brs,SarSpi2}.

\section{Multiple elliptic hypergeometric integrals}

There are many multiple integral generalizations of the elliptic
beta integral evaluation and of the $V$-function. For all of them
the integrands satisfy a set of linear $q$-difference
equations of the first order in the integration variables
with the elliptic function coefficients, similar to the univariate case.

We present the most useful examples of integrals associated with the
root systems $C_n$ and $A_n$. In \cite{die-spi:selberg} it was suggested
to distinguish two types of the multiple elliptic beta integrals: those
for which the number of parameters depends on the rank of the root system
were tagged as type I, and for type II this number is fixed. There is also
a difference in the methods of proving their evaluation formulas.

So, the type I integral on the $C_n$ root system has the following form.
Take $2n+4$ complex parameters $t_1,\ldots,t_{2n+4}$ and bases $p,q$
with the absolute values $|p|, |q|, |t_j|<1$, and impose the
balancing condition $\prod_{j=1}^{2n+4}t_j=pq$. Then one has
\begin{eqnarray}\nonumber
&& \kappa_n\int_{\T^n}\prod_{1\leq j<k\leq n}\frac{1}{\eg(z_j^{\pm 1} z_k^{\pm 1};p,q)}
\prod_{j=1}^n\frac{\prod_{i=1}^{2n+4}\eg(t_iz_j^{\pm 1};p,q)}
{\eg(z_j^{\pm2};p,q)}\prod_{j=1}^n\frac{dz_j}{z_j}
\\ && \makebox[4em]{}
=\prod_{1\leq i<j \leq 2n+4}\eg(t_it_j;p,q), \qquad
\kappa_{n}=\frac{(p;p)_\infty^n(q;q)_\infty^n}{(4\pi {i})^n n!}.
\label{C-typeI}\end{eqnarray}

The simplest proof of this relation uses a direct generalization of the method
described above for the univariate case. The ratio of the integrand and the right-hand
side expression satisfies a linear difference equation in parameters and integration
variables similar to \eqref{keyeq}. Other univariate arguments generalize as well \cite{spi:short},
which yields \eqref{C-typeI}. The original work \cite{die-spi:selberg}, where this formula
was suggested, contained only its partial justification. The first complete proof was given by
Rains \cite{rai:trans} using a different method and in a substantially more general setting.
Namely, the following transformation formula was established in \cite{rai:trans}
\begin{equation}
I_n^{(m)}(t_1,\ldots,t_{2n+2m+4})=\prod_{1\leq i<j\leq 2n+2m+4}\eg(t_it_j;p,q)\;
I_m^{(n)}\left(\frac{\sqrt{pq}}{t_1},\ldots,\frac{\sqrt{pq}}{t_{2n+2m+4}}\right)
\label{trafo}\end{equation}
for the integrals
$$
I_n^{(m)}(\underline{t})=\kappa_n\int_{\T^n}
\prod_{1\leq i<j\leq n} \frac{1}{\eg(z_i^{\pm 1}z_j^{\pm 1};p,q)}
\prod_{j=1}^n\frac{\prod_{i=1}^{2n+2m+4}\eg(t_iz_j^{\pm 1};p,q)}
{\eg(z_j^{\pm 2};p,q)}\frac{dz_j}{z_j},
$$
where $|t_j|<1$ and $\prod_{j=1}^{2n+2m+4}t_j=(pq)^{m+1}$.
As  shown in  \cite{rai:trans}, analytically the product
$\prod_{1\leq k<l\leq 2n+2m+4}(t_kt_l;p,q)_\infty
I_n^{(m)}(\underline{t})$ is a holomorphic function of its parameters.
Relation \eqref{trafo} can be considered as an elliptic analogue
of the  symmetry transformation for ordinary hypergeometric integrals
established by Dixon \cite{Dixon}.  Clearly it represents a
multivariable extension of the third $V$-function symmetry transformation  \eqref{E7-3}.

 In \cite{RS} these
integrals were written as determinants of univariate integrals
\begin{eqnarray*} &&\makebox[-1em]{}
I_n^{(m)}(t_1,\ldots,t_{2n+2m+4})= \prod_{1\leq i<j\leq n}
\frac{1}{a_j\theta(a_ia_j^{\pm 1};p)b_j\theta(b_ib_j^{\pm 1};q)}
\\ && \makebox[-1em]{} \times
\det_{1\le i,j\le n}\left(\kappa\int_\T \frac{\prod_{r=1}^{2n+2m+4}
\eg(t_r z^{\pm 1};p,q)}{\eg(z^{\pm2};p,q)}
\prod_{k\ne i} \theta(a_kz^{\pm 1};p)
\prod_{k\ne j} \theta(b_kz^{\pm 1};q)
\; \frac{dz}{z} \right),
 \end{eqnarray*}
where $a_i, b_i$ are arbitrary auxiliary variables. Curiously,  the Casoratian
\eqref{V-det} emerges here as the required determinant for the choice $n=2, m=0$
and $a_i=b_i=t_i$, which yields the evaluation formula \eqref{C-typeI} for $n=2$.

For the description of type II $C_n$ elliptic beta integral
introduced in \cite{die-spi:elliptic} one needs
seven complex parameters $t, t_a$, $a=1,\ldots , 6$, and bases $p$ and $q$
lying inside the unit disk $|p|, |q|,$ $|t|,$ $|t_a| <1,$ and
satisfying the balancing condition $t^{2n-2}\prod_{i=1}^6t_i=pq$. Then
the following integral evaluation holds true
\begin{eqnarray}\nonumber &&
\kappa_n\int_{\T^n} \prod_{1\leq j<k\leq n}
\frac{\eg(tz_j^{\pm 1} z_k^{\pm 1};p,q)}{\eg(z_j^{\pm 1} z_k^{\pm 1};p,q)}
\prod_{j=1}^n\frac{\prod_{i=1}^6\eg(t_iz_j^{\pm 1};p,q)}{\eg(z_j^{\pm2};p,q)}
\prod_{j=1}^n\frac{dz_j}{z_j}
\\ && \makebox[4em]{}
= \prod_{j=1}^n\left(\frac{\eg(t^j;p,q)}{\eg(t;p,q)}
\prod_{1\leq i<k\leq 6}\eg(t^{j-1}t_it_k;p,q )\right).
\label{C-typeII}\end{eqnarray}

As mentioned, the type II integral can be proved  by a different
method than the type I case \cite{die-spi:selberg}.
Assuming that $t_6$ is a dependent
variable, we denote the integral on the left-hand side of \eqref{C-typeII}
as $I_n(t,t_1,\ldots,t_5)$ and consider the $(2n-1)$-fold integral
\begin{eqnarray}
&& \makebox[-1em]{}
\int_{\T^{2n-1}}
\prod_{1\leq j<k\leq n}\frac{1}{\eg(z_j^{\pm 1} z_k^{\pm 1};p,q)}
 \prod_{l=1}^n\frac{\prod_{r=0}^5\eg(t_rz_l^{\pm 1};p,q)}
{\eg(z_l^{\pm2};p,q)}\frac{dz_l}{z_l}
\nonumber \\ &&
\times
\prod_{\stackrel{1\leq j\leq n}{1\leq k\leq n-1}}
\eg(t^{1/2}z_j^{\pm 1} w_k^{\pm 1};p,q)
\prod_{1\leq j<k\leq n-1}\frac{1}{\eg(w_j^{\pm 1} w_k^{\pm 1};p,q)}\nonumber \\
&&  \makebox[2em]{}
\times \prod_{j=1}^{n-1}
\frac{\eg(w_j^{\pm 1} t^{n-3/2}\prod_{s=1}^5t_s;p,q)}
{\eg(w_j^{\pm2},w_j^{\pm 1} t^{2n-3/2}\prod_{s=1}^5t_s;p,q)}
\frac{dw_j}{w_j},
\nonumber
\end{eqnarray}
where we introduced an auxiliary variable $t_0$ via the
relation $t^{n-1}\prod_{r=0}^5t_r=pq$.

Integrals over $w_j$ or $z_j$ can be computed explicitly
using the type I $C_n$-integral \eqref{C-typeII}.
Doing these integrations in different orders, one obtains the recurrence relation:
$$ \makebox[-1em]{}
I_n(t,t_1,\ldots,t_5)= \frac{\Gamma(t^n;p,q)}{\Gamma(t;p,q)}\makebox[-0.5em]{}
\prod_{0\leq r<s\leq 5}\makebox[-0.5em]{}\Gamma(t_rt_s;p,q)\;
I_{n-1}(t,t^{1/2}t_1,\ldots,t^{1/2}t_5)
$$
with known $n=1$ initial condition. Resolving this recurrence
one comes to the desired formula.

Expressing one of the parameters $t_i$ in terms of others using
the balancing condition and taking the limit $p\to 0$ for fixed
values of independent parameters, one reduces the above integrals
to Gustafson's $C_n$ $q$-beta integrals from \cite{gust2}.
Relation \eqref{C-typeII} has a meaning of an elliptic extension
of the Selberg integral evaluation formula \cite{aar,for-war},
which emerges as a result of its sequential degenerations.

Let us present also an elliptic beta integral of type I for the $A_n$
root system suggested in \cite{spi:theta2} and proven in \cite{rai:trans} and \cite{spi:short}.
Take $2n+4$ parameters $t_m, s_m,\, m=1,\ldots, n+2,$ and bases $p, q$
satisfying the constraints $|p|, |q|, |t_m|, |s_m|<1$ and the balancing
condition $ST=pq$, where $S=\prod_{m=1}^{n+2}s_m$ and $T=\prod_{m=1}^{n+2}t_m$.
Then the following integral can be computed explicitly
\begin{eqnarray}\nonumber
&& \mu_{n}\int_{{\mathbb T}^n}
\prod_{1\le j<k\le n+1}\frac{1}{\eg(z_jz_k^{-1},z_j^{-1}z_k;p,q)}
\,\prod_{j=1}^{n+1}\prod_{m=1}^{n+2}\eg(s_mz_j,t_mz_j^{-1};p,q)
\,\prod_{j=1}^n\frac{dz_j}{z_j}
\\ && \makebox[4em]{}
=\prod_{m=1}^{n+2} \eg(Ss_m^{-1},Tt_m^{-1};p,q)
\prod_{k,m=1}^{n+2} \eg(s_kt_m;p,q),
\label{AI}\end{eqnarray}
where $z_1z_2\cdots z_{n+1}=1$ and
$$
\mu_{n}=\frac{(p;p)_\infty^n(q;q)_\infty^n}{(2\pi i)^n(n+1)!}.
$$

Relations to the root systems emerge from the following observation.
Combinations of the integration variables of the form
$z_j^{\pm1}z_k^{\pm1},\, j<k$, $z_j^{\pm2}$  in \eqref{C-typeI}, \eqref{C-typeII}
and $z_jz_k^{-1},\, j\neq k$, in \eqref{AI} can be identified
with formal exponentials of the roots $\pm e_j \pm e_k,\, j<k$, $\pm2 e_j$
and $e_j-e_k,\, j\neq k,$ of the $C_n$ and $A_n$ root systems, respectively.

\section{Rarefied elliptic hypergeometric integrals}

Recently a further modification of the elliptic hypergeometric integrals
has been introduced in \cite{kels,KY,RW,spi:rare}. It emerged  from
considerations of supersymmetric quantum field theories on particular four dimensional
space-time background $S^1\times L(r,-1)_\tau$ involving a special lens space.
The general squashed lens space $L(r,k)_\tau$ is obtained from the
squashed three-dimensional sphere
in the complex representation $|\tau z_1|^2+|\tau^{-1}z_2|^2=1$
by identification of the points $(e^{2\pi i /r}z_1, e^{2\pi i k/r}z_2)\sim (z_1,z_2)$
for $k,\, r$ positive coprime integers $0<k<r$. Let us describe briefly
corresponding generalizations of the elliptic hypergeometric identities.

A proper extension of the elliptic gamma function, associated with a
special lens space, is determined
by two standard elliptic gamma functions with different bases
\begin{eqnarray}\label{reg} &&
\gamma^{(r)}(z,m;p,q):=\Gamma(zp^m;p^r, pq)\Gamma(zq^{r-m};q^r, pq),
\nonumber\end{eqnarray}
where one has two integer parameters $r\in\mathbb{Z}_{>0}$ and $m\in\mathbb{Z}$.
Using the double elliptic gamma function $\Gamma(z;p,q,t)$
with a special choice $t=pq$, one can  write
\begin{eqnarray} \nonumber &&
\gamma^{(r)}(z,m;p,q)
=\frac{\Gamma(q^rz p^m;p^r,q^r,pq)}{\Gamma(z p^m;p^r,q^r,pq)}
\frac{\Gamma(p^r z q^{r-m};p^r,q^r,pq)}{\Gamma(z q^{r-m};p^r,q^r,pq)}
\\  && \makebox[6em]{}
=\frac{\Gamma((pq)^mq^{r-m}z;p^r,q^r,pq)}{\Gamma(q^{r-m}z;p^r,q^r,pq)}
\frac{\Gamma((pq)^{r-m}p^mz;p^r,q^r,pq)}{\Gamma(p^mz;p^r,q^r,pq)},
\label{regm_egm3}\end{eqnarray}
which yields the ``rarefied" product representation for $\gamma^{(r)}(z,m;p,q)$
of the form
\beq
\gamma^{(r)}(z,m;p,q)
=\prod_{k=0}^{m-1}\Gamma(q^{r-m}z(pq)^k;p^r,q^r)\prod_{k=0}^{r-m-1}\Gamma(p^mz(pq)^k;p^r,q^r),
\label{regm_r1}\ee
valid for $0\leq m\leq r$ (similar expression exists for other values of $m$).
This function is quasiperiodic in the discrete variable
\beq
\gamma^{(r)}(z,m+r;p,q)=(-z)^{-m}q^{m(m+1)/2}p^{-m(m-1)/2}\gamma^{(r)}(z,m;p,q).
\label{r-per_gen}\ee

The normalized function
\beq
\Gamma^{(r)}(z,m;p,q)
:=\left(-\frac{z}{\sqrt{pq}}\right)^{\frac{m(m-1)}{2}}
\left(\frac{p}{q}\right)^{\frac{m(m-1)(2m-1)}{12}}\gamma^{(r)}(z,m;p,q).
\label{regm}\ee
was called in \cite{spi:rare} the rarefied elliptic gamma function.
For $r=1$ independently of $m$ one has $\Gamma^{(1)}(z,m;p,q)=\Gamma(z;p,q)$, which
provides a very convenient verification of identities involving $\Gamma^{(r)}(z,m;p,q)$.

Let us describe some properties of this function. The $p, q$ permutational symmetry
changes to
\beq
\Gamma^{(r)}(z,m;p,q)=\Gamma^{(r)}(z,-m;q,p).
\label{sympq}\ee
Instead of the plain difference equations one has simple recurrence relations
\begin{eqnarray} \nonumber &&
\Gamma^{(r)}(qz,m+1;p,q)=(-z)^m p^{\frac{m(m-1)}{2}}\theta(zp^{m};p^r)\Gamma^{(r)}(z,m;p,q),
\\ &&
\Gamma^{(r)}(pz,m-1;p,q)=(-z)^{-m}q^{\frac{m(m+1)}{2}}
\theta(zq^{-m};q^r)\Gamma^{(r)}(z,m;p,q).
\label{Gp}\end{eqnarray}
The inversion relation takes the form
\beq\label{invGamma}
\Gamma^{(r)}(z,m;p,q)\Gamma^{(r)}(\textstyle{\frac{pq}{z}},-m;p,q)=1,
\ee
and the limiting relation needed for the residue calculus reads
 \beq\label{rfres} \qquad
\stackreb{\lim}{z\to 1}(1-z) \Gamma^{(r)}(z,0;p,q)
=\stackreb{\lim}{z\to 1}(1-z) \gamma^{(r)}(z,0;p,q)
=\frac{1}{(p^r;p^r)_\infty(q^r;q^r)_\infty}.
 \ee

{\em Exercise:} verify all these relations.

The rarefied version of the elliptic beta integral has the following form.
We take continuous parameters $t_1, \dots ,t_6,p,q $ and discrete ones $n_1,\ldots, n_6\in\mathbb{Z}+\nu$,
where $\nu=0,\frac{1}{2}$,
satisfying the constraints $|t_a|,|p|,|q| <1$ and the balancing condition
$$
\prod_{a=1}^6 t_a=pq, \qquad \sum_{a=1}^6 n_a=0.
$$
Then
\begin{equation} \label{rfint}
\kappa^{(r)}\sum_{m\in\Z_r+\nu}
\int_{\mathbb{T}}\rho^{(r)}(z,m;\underline{t},\underline{n})
\frac{dz}{z}
= \prod_{1 \leq a < b \leq 6} \Gamma^{(r)}(t_a t_b,n_a+n_b;p,q),
\end{equation}
where $\mathbb T$ is the positively oriented unit circle,
$$
\kappa^{(r)}=\frac{(p^r;p^r)_\infty (q^r;q^r)_\infty}{4\pi {i}},
$$
and the integrand has the form
\beq
\rho^{(r)}(z,m;\underline{t},\underline{n}):=
\frac{\prod_{a=1}^6\Gamma^{(r)}(t_az^{\pm1},n_a \pm m;p,q)}
{\Gamma^{(r)}(z^{\pm2},\pm 2m);p,q)}.
\label{ker_rebeta}\ee
Here we use the compact notation
\beq
\Gamma^{(r)}(tz^{\pm1},n\pm m;p,q):=\Gamma^{(r)}(tz,n+m;p,q)
\Gamma^{(r)}(tz^{-1},n-m;p,q).
\ee

For $r=1$ one gets relation \eqref{elbeta} and the $r>1,\,\nu=0$  case
of the evaluation \eqref{rfint} was established by Kels in \cite{kels},
for $r>1,\,\nu=\frac{1}{2}$ it was proven in \cite{spi:rare} in the presented form
and in \cite{KY} in the equivalent form of $A_1\leftrightarrow A_0$ symmetry
transformation.

A good calculational exercise is the proof of periodicity
\beq\label{period}
\rho^{(r)}(z,m+r;\underline{t},\underline{n})=\rho^{(r)}(z,m;\underline{t},\underline{n}),
\ee
because of which the sum over $m-\nu=0,1,\ldots, r-1$ is equal to sums over any $r$
consecutive values of $m$. There is a particular symmetry between
the terms in this sum following from the obvious relation
$$
\rho^{(r)}(z,-m;\underline{t},\underline{n})=\rho^{(r)}(z^{-1},m;\underline{t},\underline{n}).
$$
Due to the $r$-periodicity in $m$ one has
\begin{eqnarray*} &&
c_{r-m}:=\int_{\mathbb{T}}\rho^{(r)}(z,r-m;\underline{t},\underline{n})\frac{dz}{z}=
\int_{\mathbb{T}}\rho^{(r)}(z^{-1},m;\underline{t},\underline{n})\frac{dz}{z}
\\ && \makebox[2,5em]{}
= \int_{\mathbb{T}}\rho^{(r)}(z,m;\underline{t},\underline{n})\frac{dz}{z}
=c_{m}.
\end{eqnarray*}
As a result the sum over $m$ in \eqref{rfint} can be written for $\nu=0$ as
\beq
\sum_{m=0}^{r-1}c_m =
\begin{cases}
c_0+c_{r/2}+2\sum_{m=1}^{r/2-1}c_m & \text{for even}\, r, \\
c_0+2\sum_{m=1}^{(r-1)/2}c_m  & \text{for odd}\, r,
\end{cases}
\label{e=0red}\ee
and for $\nu=\frac{1}{2}$ as
\beq
\sum_{m=1/2}^{r-1/2}c_m =
\begin{cases}
2\sum_{m=1/2}^{r/2-1/2}c_m & \text{for even}\, r, \\
c_{r/2}+2\sum_{m=1/2}^{(r-2)/2}c_m  & \text{for odd}\, r.
\end{cases}
\label{e=1red}\ee

The type I multiple rarefied elliptic beta integral for the root system $C_n$
has the form
\begin{eqnarray}  && \makebox[-2em]{}
\kappa_{n}^{(r)} \sum_{m_1,\ldots,m_n\in \mathbb{Z}_k+\nu}
\int_{\mathbb{T}^n}\rho_{I}(z_j,m_j;\underline{t},\underline{n})
\prod_{j=1}^n\frac{dz_j}{z_j}
= \prod_{1 \leq a < b \leq 2 n+4} \Gamma(t_a t_b,n_a+n_b;p,q),
\label{rfintCI}\end{eqnarray}
where $\mathbb T$ is the unit circle of positive orientation, $\kappa_{n}^{(r)}$ is obtained from $\kappa_n$
 after replacing $p, q \to p^r, q^r$, and the kernel is
\begin{eqnarray} \nonumber &&
\rho_{I}({z}_j,{m}_j;\underline{t},\underline{n}):=\prod_{1\leq j<k\leq n}
\frac{1}{\Gamma(z_j^{\pm 1}z_k^{\pm 1},\pm m_j\pm m_k)}
\prod_{j=1}^n\frac{\prod_{a=1}^{2n+4}
\Gamma(t_a^{\pm1} z_j,n_a\pm m_j)}
{\Gamma(z_j^{\pm 2},\pm 2m_j)},
\label{ker_ebetaCI}\end{eqnarray}
where parameters $ t_a, z_j\in \mathbb{C}^\times,\, n_a, m_j\in\mathbb{Z}+\nu,$
satisfy the constraints $|t_a|<1$ and the balancing condition
\beq
\prod_{a=1}^{2n+4}t_a=pq,\qquad \sum_{a=1}^{2n+4}n_a=0.
\label{balI}\ee

The proof of the univariate case $n=1$ can be adapted to the present situation
by adjoining the peculiarities characteristic to the proof of type I integral \eqref{C-typeI}
as well as the $r$-periodicity of the kernel in the discrete summation variables.

Similarly one can construct a computable rarefied type II $C_d$-integral,
where for convenience we denoted the rank of the root system as $d$.
For that it is necessary to take continuous parameters $t, t_a \in \mathbb{C}^\times, a=1,\ldots , 6,$
and bases $p, q$ such that  $|p|, |q|,$ $|t|,$ $|t_a| <1$. Additionally, one needs
eight discrete variables $n\in\Z, n_a\in\mathbb{Z}+\nu,$ all together satisfying
the balancing condition
\begin{equation}
t^{2d-2}\prod_{a=1}^6t_a=pq, \qquad 2n (d-1)+\sum_{a=1}^6n_a=0.
\label{balII}\end{equation}
 Then
\begin{eqnarray}\nonumber &&\makebox[0em]{}
\kappa_d^{(r)}\sum_{m_1,\ldots,m_d\in\Z_r+\nu}\int_{\T^d} \prod_{1\leq j<k\leq d}
\frac{\Gamma (tz_j^{\pm 1} z_k^{\pm 1},n\pm m_j\pm m_k)}
{\Gamma (z_j^{\pm 1} z_k^{\pm 1},\pm m_j\pm m_k)}
\prod_{j=1}^d\frac{\prod_{a=1}^6\Gamma(t_az_j^{\pm 1},n_a\pm m_j)}
{\Gamma (z_j^{\pm2},\pm 2m_j)}
\frac{dz_j}{z_j}
\\  && \makebox[4em]{}
= \prod_{j=1}^d\left(\frac{\Gamma (t^j, nj)}{\Gamma (t,n)}
\prod_{1\leq a<b\leq 6}\Gamma (t^{j-1}t_at_b, n(j-1)+n_a+n_b)\right).
\label{rfintCII}\end{eqnarray}
This formula is proved in a way similar to the $r=1$ case \eqref{C-typeI}, i.e.
by considering a $(2d-1)$-fold combination of summations and integrations of
a specific function admitting usage of the rarefied  type I $C_d$-formula \eqref{rfintCI}
in two different sets of discrete summation and continuous integration variables
which establishes a recurrence relation in the rank of the root system. For more details
on these results, as well as generalizations of the $V$-function and
elliptic hypergeometric equation, see \cite{spi:rare}. Symmetry
transformations for some multidimensional elliptic hypergeometric integrals
were extended to the rarefied case in \cite{KY}.

The rarefied hyperbolic hypergeometric integrals for general lens space were
discussed in \cite{dimofte,Gor}. In particular, in \cite{Gor} a general univariate
computable rarefied hyperbolic beta integral evaluation formula has been established.

\section{An integral Bailey lemma}

Using properties of the elliptic beta integral, the following
integral transformation was introduced in \cite{spi:bailey2}
\begin{equation}
\beta(w,t)=M(t)_{wz}\alpha(z,t):=\frac{(p;p)_\infty(q;q)_\infty}{4\pi{i}}\int_\mathbb{T}
\frac{\Gamma(tw^{\pm1}z^{\pm1};p,q)}
{\Gamma(t^2,z^{\pm2};p,q)}\alpha(z,t)\frac{dz}{z}
\label{EFT}\end{equation}
with the assumption that $|tw^{\pm1}|<1$. The latter constraints can be relaxed by
analytic continuation, e.g. by deforming the contour of integration, provided no
singularities of the integrand are crossed during such a deformation.
Pairs of functions connected by \eqref{EFT} were called integral elliptic Bailey pairs with respect
to the parameter $t$. Using the evaluation formula \eqref{elbeta} one can find a particular
explicit Bailey pair $\alpha(z,t)$ and $\beta(z,t)$.
Such a terminology emerged from the theory of Bailey chains
providing a systematic tool for constructing nontrivial identities for $q$-hypergeometric
series \cite{Warnaar}. In particular, it was targeted at the proof of Rogers-Ramanujan type
identities. The definition \eqref{EFT} yielded the very first generalization
of the Bailey chains technique from series to integrals.

As shown in \cite{spi-war:inversions}, on the space of $A_1$-symmetric functions
 $f(z)=f(z^{-1})$ under particular constraints on the parameters
and appropriate choice of the integration contours for analytically continued
operators, the operators $M(t^{-1})_{wz}$ and $M(t)_{wz}$
become inversions one of the other.  Passing to the real line integration
one can use the generalized functions
and symbolically write $M(t^{-1})M(t)=1$,  where 1 means
an integral operator with the Dirac delta-function kernel \cite{S5,stat}.
It is due to this $t\to t^{-1}$ inversion relation, which looks
similar to the Fourier transform, that the transformation \eqref{EFT}
is referred to as the ``elliptic Fourier transformation". Another similarity is that
in both cases some nontrivial operators --- the derivative and $q$-scaling are
converted to the multiplication by a function --- the linear and theta functions, respectively.

Let us indicate how the true Fourier transform actually emerges in a particular degeneration limit
of \eqref{EFT}. Take first the limit $p\to 0$ for fixed $q,t$ and $w$. This yields
$$
\beta(w,t)=\frac{(q;q)_\infty}{4\pi{i}}\int_\mathbb{T}
\frac{(t^2,z^{\pm2};q)_\infty}{(tw^{\pm1}z^{\pm1};q)_\infty}
\alpha(z,t)\frac{dz}{z}.
$$
In the integrand one can write
$$
(z^{\pm2};q)_\infty=(z^{\pm1},-z^{\pm1};q)_\infty (q^{1/2}z^{\pm1},-q^{1/2}z^{\pm1};q)_\infty.
$$
We can rewrite the above transform in a renormalized form
$$
\tilde\beta(w)=\frac{1}{4\pi{i}}\int_\mathbb{T}
\frac{(z^{\pm1},-z^{\pm1};q)_\infty}{(tw^{\pm1}z^{\pm1};q)_\infty}\tilde\alpha(z)\frac{dz}{z},
$$
where $\tilde\beta(w):=\beta(w,t)/(q,t^2;q)_\infty$ and $\tilde\alpha(z):=(q^{1/2}z^{\pm1},-q^{1/2}z^{\pm1};q)_\infty\alpha(z,t)$.
Passing to the angular parameter $\theta$, $z=e^{i\theta}$, and introducing a new integration
variable $x=\cos \theta$ we obtain $\int_\T dz/z=2i\int_{-1}^1dx/\sqrt{1-x^2}$.
Denoting $tw=q^{\alpha+1/2}, tw^{-1}=-q^{\beta+1/2}$,
and using the $q$-binomial limiting relation \eqref{limbin}
we deduce
$$
\stackreb{\lim}{q\to1^-} \frac{(z^{\pm1},-z^{\pm1};q)_\infty}{(q^{\alpha+1/2}z^{\pm1},
-q^{\beta+1/2}z^{\pm1};q)_\infty}=2^{\alpha+\beta+1}(1-x)^{\alpha+1/2}(1+x)^{\beta+1/2}
$$
and come to the integral transform
$$
g(\alpha,\beta)=\frac{2^{\alpha+\beta}}{\pi}\int_{-1}^1(1-x)^\alpha(1+x)^\beta f(x)dx.
$$
where we have to assume that Re$(\alpha)$, Re$(\beta)>-1$ for convergence
of the integral for regular functions $f(x)$.
Rescaling in the integral $x\to x/\lambda$ and taking the limit $\lambda\to+\infty$,
we obtain asymptotically the transform
$$
\frac{2^{\alpha+\beta}}{\pi\lambda}\int_{-\infty}^\infty e^{\frac{\beta-\alpha}{\lambda} x -\frac{\alpha+\beta}{2\lambda^2} x^2 +O(\frac{\beta-\alpha}{\lambda^3})} \tilde f(x)dx,
$$
where $\tilde f(x)=\lim_{\lambda\to\infty}f(x/\lambda)$.
Demanding that $\alpha+\beta=o(\lambda^2)$ and $\beta-\alpha=iy\lambda$ for some finite variable $y$,
we obtain
$$
\frac{2^{\alpha+\beta}}{\pi\lambda}\int_{-\infty}^\infty e^{iyx} \tilde f(x)dx,
$$
which is the standard Fourier transformation up to some diverging factor. So, in terms of the
original variables, the action of the integral operator \eqref{EFT} passes to the ordinary Fourier transformation
after setting $p=0$, proper normalization of the source and image functions,
and taking the limit $q\to 1^-$ in the parameterization $w=-iq^{-iy\lambda/2}$,
$z+z^{-1}=2x/\lambda$, $t=iq^{c}$ with the subsequent limit $\lambda\to+\infty$ and
the constraint that $c=(\alpha+\beta+1)/2$ is an arbitrary
parameter which may grow only slower than $\lambda^2$.

The integral Bailey lemma provides an
algorithm for constructing infinitely many Bailey pairs from a given one.
It is formulated as follows. Let $\alpha(z,t)$ and $\beta(z,t)$ be some
functions related by \eqref{EFT} for some parameter $t$. Then the functions
\begin{eqnarray} \nonumber &&
\alpha'(w,st)=D(s;y,w)\alpha(w,t),\quad D(s;y,w)
=\Gamma(\sqrt{pq}s^{-1}y^{\pm1}w^{\pm1};p,q),
\\  &&
\beta'(w,st)=D(t^{-1};y,w) M(s)_{wx}D(st;y,x)\beta(x,t),
\label{BL}\end{eqnarray}
where $w\in \T$, $|s|,|t|<1, |\sqrt{pq}y^{\pm1}|<|st|$, form an integral
elliptic Bailey pair with respect to the parameter $st$. Note that the
parameters $s$ and $y$ are two new arbitrary variables.

It is necessary to show that $\beta'(w,st)=M(st)_{wz}\alpha'(z,st)$.
Substitute in both sides of this equality the definitions \eqref{BL}
and use the relation $D(t^{-1};y, w)=1/D(t;y,w)$ following from
the elliptic gamma function inversion property. This yields the
operator identity
\begin{equation}
M(s)_{wx}D(st;y,x)M(t)_{xz}=D(t;y,w)M(st)_{wz}D(s;y,z).
\label{STR}\end{equation}
Substitution of the explicit forms of $M$- and $D$-operators shows that
the integral over the variable $x$ on the left-hand side of \eqref{STR} can
be computed explicitly using the elliptic beta integral evaluation formula.
The resulting expression takes exactly the form given on the right-hand side.

Iterative applications of the maps \eqref{BL} lead to a chain of Bailey pairs
satisfying by definition the key relation \eqref{EFT}. Explicitly this leads
to certain nontrivial identities for multiple elliptic hy\-per\-geo\-met\-ric
integrals. For instance, if the pair $\alpha$ and $\beta$ is determined
from the formula \eqref{elbeta}, then the relation $\beta'(w,st)=M(st)_{wz}\alpha'(z,st)$
yields the key $W(E_7)$-transformation for the $V$-function \eqref{E7-1}.

As shown in \cite{DS} the algebraic relations emerging from the described integral
Bailey lemma can have the meaning of Coxeter relations for a permutation group.
For that interpretation we introduce three operators $\mathrm{S}_{1,2,3}(\mathbf{t})$
acting on the functions of two complex variables $f(z_1,z_2)$ as follows
\begin{eqnarray*} && 
[\mathrm{S}_1(\mathbf{t})f](z_1,z_2):=M(t_1/t_2)_{z_1z}f(z,z_2), \quad
\\ &&
[\mathrm{S}_2(\mathbf{t})f](z_1,z_2):=D(t_2/t_3;z_1,z_2)f(z_1,z_2),
\\ &&
[\mathrm{S}_3(\mathbf{t})f](z_1,z_2):=M(t_3/t_4)_{z_2z}f(z_1,z),
\end{eqnarray*}
for some complex parameters $\mathbf{t}=(t_1,t_2,t_3,t_4)$.
The products of these operators are defined via the cocycle condition
$$
\mathrm{S}_j\mathrm{S}_k:=\mathrm{S}_j(s_k(\mathbf{t}))\mathrm{S}_k(\mathbf{t}),
$$
where $s_k$ are elementary transposition operators generating
the permutation group $\mathfrak{S}_4$:
$$
s_1(\mathbf{t})=(t_2,t_1,t_3,t_4), \quad
s_2(\mathbf{t})=(t_1,t_3,t_2,t_4), \quad
s_3(\mathbf{t})=(t_1,t_2,t_4,t_3).
$$
Now one can check validity of the Coxeter relations
\begin{equation}
\mathrm{S}_j^2=1, \quad \mathrm{S}_i\mathrm{S}_j=\mathrm{S}_j\mathrm{S}_i \
\text{ for } \ |i-j|>1, \quad
\mathrm{S}_j\mathrm{S}_{j+1}\mathrm{S}_j
=\mathrm{S}_{j+1}\mathrm{S}_j\mathrm{S}_{j+1}
\label{coxeter}\end{equation}
as a consequence of properties of the Bailey lemma operator entries.
The quadratic relations represent inversion relations for
the $M$- and $D$-operators. The cubic relation is equivalent to \eqref{STR}
and it is called also the star-triangle relation. A somewhat different application
of the operator identity \eqref{STR} is considered in \cite{rai:sigma}.
Extension of the above considerations to the rarefied elliptic beta integral was
considered in \cite{spi:rareYBE}.

Let us replace in \eqref{STR} all variables $z\to e^{ i z}$,
$x\to e^{ i x}$, $y\to e^{ i y}$, $w\to e^{ i w}$ and denote
$s=e^{-\alpha}$, $t=e^{-\beta}$, $pq= e^{-2\eta}$,
and pass to the integrations over the line segment $x,z\in [0,2\pi]$.
Applying now this operator identity to the Dirac delta-function $(\delta(z-u)+\delta(z+u))/2$
for some parameter $u$, one comes to formula \eqref{elbeta} written in the form
\begin{eqnarray}\nonumber
&& \int_{0}^{2\pi} \rho(x)D_{\eta-\alpha}(w,x)D_{\alpha+\beta}(y,x)D_{\eta-\beta}(u,x)dx
\\ && \makebox[2em]{}
=\chi(\alpha, \beta)D_\beta(y,w)D_{\eta-\alpha-\beta}(w,u)D_{\alpha}(y,u),
\label{astr}\end{eqnarray}
where
\begin{equation}
D_{\alpha}(y,u)=D(e^{-\alpha}; e^{ {i}  y},
e^{{i} u})=\Gamma(e^{ \alpha-\eta\pm iy \pm iu)};p,q)
\label{weight}\end{equation}
and
\begin{eqnarray*}
&&
\rho(u)=\frac{ (p;p)_\infty(q;q)_\infty}{2\Gamma(e^{\pm 2i u};p,q)},
\quad \chi(\alpha,\beta)= \Gamma(e^{-2\alpha},e^{-2\beta},e^{2\alpha+2\beta-2\eta};p,q).
 \end{eqnarray*}

In \cite{BS} this form of the star-triangle relation was
used for building a new two-dimensional integrable lattice model. Namely, one considers
a two-dimensional square lattice and ascribes the Boltzmann weight $D_\alpha(x,u)$
to the horizontal edges connecting continuous spins $x$ and $u$ sitting in the
neighboring vertices of the lattice. The vertical edges have Boltzmann weights
$D_{\eta -\alpha}(x,u)$. Each vertex has the self-interaction energy $\rho(u)$.

Let us substitute in \eqref{astr} $D_{\alpha}(y,w)=m(\alpha)W_{\alpha}(y,w)$ and choose the
normalization constant $m(\alpha)$ from the condition
\beq
\frac{m(\alpha)m(\beta)m(\eta-\alpha-\beta)}
{m(\eta-\alpha)m(\eta-\beta)m(\alpha+\beta)}
\chi(\alpha, \beta)=1.
\lab{chi}\ee
This gives a compact block representation of the elliptic beta integral evaluation
$$
\int_{0}^{2\pi} \rho(x)W_{\eta-\alpha}(w,x)W_{\alpha+\beta}(y,x)W_{\eta-\beta}(u,x)dx
=W_\beta(y,w)W_{\eta-\alpha-\beta}(w,u)W_{\alpha}(y,u).
$$
Equality \eqref{chi} holds true, if
$$
m(\alpha+\eta)=\Gamma(e^{2\alpha};p,q)m(-\alpha).
$$

In order to compute $m(\alpha)$ it is convenient to consider the function
\begin{equation}
\mu(x;p,q,t)
=\frac{\Gamma(xt\sqrt{pqt};p,q,t^2)}{\Gamma(x^{-1}t\sqrt{pqt};p,q,t^2)}
=\exp\Big(\sum_{n\in\Z/\{0\}} \frac{(\sqrt{pqt}x)^n}{n(1-p^n)(1-q^n)(1+t^n)}\Big),
\label{normell}\end{equation}
where $\Gamma(z;p,q,t^2)$ is the  second order elliptic gamma function with bases $p, q, t^2$.
One has the evident reflection equation $\mu(x^{-1};p,q,t)\mu(x;p,q,t)=1.$
Another easily verifiable functional equation,
$$
\mu(x;p,q,t)\mu(t^{-1}x;p,q,t)=\Gamma\Big(x\sqrt{\frac{pq}{t}};p,q\Big),
$$
becomes equivalent to the equation for $m(\alpha)$ after setting $t=pq$ and
denoting $x=e^{2\alpha}$. As a result, we find the normalizing factor of interest
\begin{equation}
m(\alpha)=\frac{\Gamma(e^{2\alpha}(pq)^2;p,q,(pq)^2)}{\Gamma(e^{-2\alpha}(pq)^2;p,q,(pq)^2)},
\quad m(\alpha)m(-\alpha)=1.
\label{ma1}\end{equation}

The partition function of the described lattice model has the form
$$
Z=\int \prod_{(ij)}W_\alpha(u_i,u_j)\prod_{(kl)}W_{\eta-\alpha}(u_k,u_l)
\prod_{m}\rho(u_m)du_m,
$$
where the product $\prod_{(ij)}$ is taken over the horizontal edges,
the product $\prod_{(kl)}$ takes into account vertical edges,
and the product in $m$ counts self-energies of all lattice vertices.
As argued in \cite{BS}, for the edge Boltzmann weights
$W_\alpha(x,u)$ the free energy per edge vanishes in the thermodynamic limit,
i.e. $\lim_{N,M\to \infty} \frac{1}{NM}\log Z=0$, where $N$ and $M$ are the numbers of
edges in the rows and columns of the lattice. As observed in \cite{stat},
the partition function $Z$ and similar ones describe superconformal indices
of four dimensional supersymmetric quiver gauge theories and the integrability conditions
represent certain electromagnetic dualities of such theories (see the next section).

The star-triangle relation can be used for constructing $R$-matrices satisfying the
 Yang-Baxter equation. We skip consideration of this subject, limiting to the
statement that the elliptic Fourier transformation operator serves as the
intertwining operator of equivalent representations of the Sklyanin algebra \cite{skl},
emerging from the $RLL$-relation associated with Baxter's 8-vertex model \cite{bax}.
More precisely, the Sklyanin algebra is generated by four operators $\mathbf{S}^a$
satisfying quadratic relations
\begin{eqnarray}\nonumber &&
\mathbf{S}^\alpha\,\mathbf{S}^\beta - \mathbf{S}^\beta\,\mathbf{S}^\alpha =
{i}\left(\mathbf{S}^0\,\mathbf{S}^\gamma +\mathbf{S}^\gamma\,\mathbf{S}^0\right),
\\ &&
\mathbf{S}^0\,\mathbf{S}^\alpha - \mathbf{S}^\alpha\,\mathbf{S}^0 =
{i}\,{J}_{\beta\gamma}\left(\mathbf{S}^\beta\,\mathbf{S}^\gamma +\mathbf{S}^\gamma\,\mathbf{S}^\beta\right),
\label{sklalg}\end{eqnarray}
where the structure constants ${J}_{\beta\gamma}=({J}_{\gamma}- {J}_{\beta})/{J}_{\alpha}$
and  $(\alpha,\beta,\gamma)$ is an arbitrary cyclic permutation of $(1,2,3)$.
An explicit realization of $\mathbf{S}^a(g)$ by finite-difference operators has been
found in \cite{skl}
\begin{eqnarray}\nonumber &&\makebox[-3em]{}
\mathbf{S}_z^a(g)= e^{\pi{i}z^2/\eta}\frac{{i}^{\delta_{a,2}}
\theta_{a+1}(\eta|\tau)}{\theta_1(2 z|\tau) } \Bigl[\,\theta_{a+1} \left(2
z-g +\eta|\tau\right)e^{\eta\partial_z}
\\  && \makebox[0em]{}
- \theta_{a+1}\left(-2z-g+\eta|\tau\right)e^{-\eta\partial_z}\Bigl]e^{-\pi{i}z^2/\eta},
\label{Sklyan} \end{eqnarray}
where $e^{\pm\eta\partial_z}$ denote the shift operators,
$e^{\pm\eta\partial_z}f(z)=f(z\pm\eta)$, and the standard theta functions are
\begin{eqnarray}
\theta_{2}(z|\tau)=\theta_1(z+{\textstyle\frac{1}{2}}|\tau), \quad
\theta_{3}(z|\tau)=e^{\frac{\pi {i}\tau}{4}+\pi {i} z}
\theta_2(z+{\textstyle \frac{\tau}{2}}|\tau), \quad
\theta_4(z|\tau)= \theta_3(z+{\textstyle\frac{1}{2}}|\tau).
\nonumber\end{eqnarray}
We added the subindex $z$ to the operators  $\mathbf{S}_z^a(g)$ in order
to indicate the arguments of the functions which they are acting on.
The usual  notation for the variable $g$ is $g=\eta(2\ell+1)$, where  $\ell\in\mathbb{C}$
is called the spin. The Casimir operators have the form
$$
\mathbf{K}_0 = \sum_{a=0}^3 \mathbf{S}^a\,\mathbf{S}^a= 4\theta_1^2\bigl(g|\tau\bigr),\quad
\mathbf{K}_2 = \sum_{\alpha=1}^3 {J}_\alpha\mathbf{S}^\alpha\,\mathbf{S}^\alpha
=4\theta_1\bigl(g-\eta|\tau\bigr)\theta_1(g+\eta|\tau).
$$
They are invariant with respect to the transformation $g\to -g$, i.e. parameters
$g$ and $-g$ correspond to equivalent representations of the Sklyanin algebra.

In \eqref{Sklyan} the operators $\mathbf{S}_z^a$ found in \cite{skl} are conjugated
by exponentials $e^{\pm \pi{i}z^2/\eta}$, which is done
for a special reason. Let us denote $q=e^{4\pi i\eta}$, $p=e^{2\pi i \tau}$, and $t= e^{-2\pi i g}$.
Then one has the following intertwining relations \cite{DS}:
\begin{equation}
{M}(t)_{WZ}\mathbf{S}_z^a(g) =
\mathbf{S}_w^a(-g){M}(t)_{WZ}, \qquad
{M}(t)_{WZ}\,\mathbf{\tilde S}_z^a(g) =
\mathbf{\tilde S}_w^a(-g)\,{M}(t)_{WZ},
\label{inter1}\end{equation}
where $W=e^{2\pi i w}$ and $Z= e^{2\pi i z}$.
The operator $M(t)_{WZ}$ is symmetric in $p$ and $q$, and the second relation in \eqref{inter1}
emerges from the first one
after interchanging  $p$ and $q$. Operators $\mathbf{\tilde S}^a_z(g) $ are thus
obtained from \eqref{Sklyan} after permutation of $2\eta$ and $\tau$ and they
realize another Sklyanin algebra with different structure constants $\tilde J_\alpha$.
Jointly these two Sklyanin algebras form the elliptic modular double \cite{S5} generalizing
Faddeev's modular double for $\mathfrak{sl}_q(2)$ algebra \cite{fad:mod}. Intertwining operators
of equivalent representations play an important role in the representation theory.
In particular, their null spaces are invariant under the action of algebra
generators which is helpful for building finite-dimensional irreducible representations.

There are useful recurrence relations for the elliptic Fourier transform operator $M(t)$
\cite{CDKK,DS2}:
\begin{equation}
{A}_k(g)\,{M}(t) = {M}(q^{-1/2}t)\,\theta_k
 \left(z | {\textstyle\frac{\tau}{2}}\right), \ \ \
{B}_k(g)\,{M}(t) = {M}\left(p^{-1/2}t\right)
\theta_k \left(z | \eta\right),
\label{RR}\end{equation}
where $k=3,4$ and ${A}_k(g)$ and ${B}_k(g)$ are the following difference operators
$$
{A}_k(g) = \frac{e^{\pi {i}\frac{(z+\eta)^2}{ \eta}}}{\theta(e^{4\pi iz};p)}
\left[ \theta_k \left(z+g+\eta | {\textstyle\frac{\tau}{2}}\right)\, e^{\eta \partial_z} -
\theta_k \left(z-g-\eta | {\textstyle\frac{\tau}{2}}\right)\, e^{-\eta \partial_z}
\right]\, e^{-\pi {i}\frac{z^2}{ \eta}},
$$
$$
{B}_k(g) =\frac{ e^{2\pi {i}\frac{(z+\tau/2)^2}{ \tau}}}
{\theta(e^{4\pi iz};q)}
\left[ \theta_k \left(z+g+{\textstyle\frac{\tau}{2}}| \eta\right)
\, e^{{\textstyle\frac{\tau}{2}} \partial_z} -
\theta_k \left(z-g-{\textstyle\frac{\tau}{2}}| \eta\right)
\, e^{-{\textstyle\frac{\tau}{2}} \partial_z}
\right]\, e^{-2\pi {i}\frac{z^2}{ \tau}}.
$$
In \eqref{RR} we drop coordinate subindices and use the convention that
the $z$-coordinate to the right of $M$-operator is the internal integration variable,
but to the left -- it is a free variable playing the role of $w$ in \eqref{inter1}.

The initial condition ${M}(1) = 1$ (the unit operator) is proved by the residue
calculus that we used in the proof of the elliptic beta integral (in this case
two pairs of poles pinch the integration contour for $t\to 1$).  Then for $t=q^{-n/2}p^{-m/2}$,
$n, m\in \Z_{\geq 0}$,
the recurrence relations can be resolved to yield the finite difference operator
\begin{eqnarray}\nonumber  &&
{M}\left(q^{-n/2}p^{-m/2}\right) =
{A}_k(n\eta-\eta+ m {\textstyle\frac{\tau}{2}})\cdots
{A}_k(\eta+m {\textstyle\frac{\tau}{2}})
{A}_k(m {\textstyle\frac{\tau}{2}})
\\  && \makebox[2em]{}
\times
{B}_k\left(m {\textstyle\frac{\tau}{2}} -{\textstyle\frac{\tau}{2}}\right)\cdots {B}_k\left({\textstyle\frac{\tau}{2}}\right)
 {B}_k(0)
\theta_k^{-m} \left(z | \eta\right)
\theta_k^{-n} \left(z | {\textstyle\frac{\tau}{2}}\right),
\label{genform}\end{eqnarray}
which does not depend on the choice of $k=3$ or $4$.
This is only one of many possible ways to represent
$M\left(q^{-n/2}p^{-m/2}\right)$ as a product of
${A}_k$- and ${B}_k$-operators.

Finally, we describe the Bailey lemma for $A_n$-root system. Define
\begin{align}\label{mop}
&M(t)_{wz}f(z):=\mu_n \int_{\mathbb{T}^n}
\frac{\prod_{j,k=1}^{n+1}\Gamma(tw_jz_k^{-1})f(z)}
{\Gamma(t^{n+1})\prod_{1\leq j<k\leq n+1}\Gamma(z_jz^{-1}_k,
z_j^{-1}z_k)}\prod_{k=1}^{n}\frac{dz_k}{2\pi{i}z_k},
\end{align}
where $\prod_{k=1}^{n+1}z_k=1$, $\Gamma(z):=\Gamma(z;p,q)$, and set
\begin{equation}
D(t;u,z):=\prod_{j=1}^{n+1}\Gamma(\sqrt{pq}t^{-\frac{n+1}{2}}\frac{u}{z_j},
\sqrt{pq}t^{-\frac{n+1}{2}}\frac{z_j}{u}), \quad D(t;u,z) D(t^{-1};u,z)=1.
\end{equation}
For $n=1$ operator \eqref{mop} coincides with \eqref{EFT}. For arbitrary $n$ it was defined in \cite{spi-war:inversions},
where the Fourier type inversion relation $M(t)_{wz}^{-1}=M(t^{-1})_{wz}$ was
established for the space of $A_n$-invariant functions under certain constraints on $t$ and $w_j$.

Similar to the univariate case, from a given Bailey pair satisfying
$\beta(w,t)=M(t)_{wz}\alpha(z,t)$, the rules
\begin{align}\nonumber
&\alpha'(w,st)=D(s;t^{-\frac{n-1}{2}}u,w)\alpha(w,t), \quad
\\
&\beta'(w,st)=D(t^{-1};s^{\frac{n-1}{2}}u,w)M(s)_{wz}D(ts;u,z)\beta(z,t)
\nonumber\end{align}
determine a new Bailey pair with respect to the parameter $st$.
From these expressions, the relation $\beta'(w,st)=M(st)_{wz}\alpha'(z,st)$ yields the cubic relation
\cite{bru-spi}
\begin{align}\label{startriangle}
&M(s)_{wz}D(st;u,z)M(t)_{zx}=
D(t;s^{\frac{n-1}{2}}u,w)M(st)_{wx}D(s;t^{-\frac{n-1}{2}}u,x),
\end{align}
which holds true due to the elliptic beta integral on the $A_n$ root system \eqref{AI}.
Although the change of $t\to  t^{-1}$ inverts $D$ and $M$ operators, for $n>1$ it is not
possible to give to equality \eqref{startriangle} a straightforward meaning of the Coxeter relation.
A substatially more complicated Bailey lemma based on the multiple $C_n$-elliptic
hypergeometric integrals of type II was formulated by Rains in \cite{rai:sigma}.

\section{Connection with four dimensional superconformal indices}

A completely unexpected development of the theory of elliptic hypergeometric
integrals emerged from quantum field theory when Dolan and Osborn \cite{DO}
have discovered that superconformal indices of four dimensional supersymmetric
gauge field theories are expressed in terms of such integrals. This was
both the most striking physical application of these integrals and
a powerful boost in understanding of their structure. We describe briefly
some ingredients of the corresponding construction and refer to surveys
\cite{Kim,RR2016} for a more detailed account and list of references.

Massless $\mathcal{N}=1$ supersymmetric field theories on the flat four dimensional space-time
have a very large symmetry group $G_{\text{full}}=SU(2,2|1)\times G\times F$.
The superconformal group $SU(2,2|1)$ contains Lorentz rotations described by
$SO(3,1)$-subgroup which is generated by $J_k, \overline{J}_k$, $k=1,2,3$.
It involves also ordinary translations and their superspace partners
generated by $P_\mu$, $\mu=0,\ldots,3$, and $Q_{\alpha}, \overline{Q}_{\dot\alpha}$,
$\alpha, \dot\alpha=1,2$, respectively. Further it includes the special conformal
transformations generator $K_\mu$ and its superpartners $S_{\alpha},\overline{S}_{\dot\alpha}$.
Finally it contains the dilations generated by $H$, and $U(1)_R$-rotations of superpartners
generated by the $R$-charge. Other symmetry groups are the
local gauge invariance group $G$ and the flavor group $F$ describing global
gauge invariance symmetries of matter superfields. Altogether they satisfy
a system of supercommutation relations forming a specific super-Lie algebra.

The superconformal index is constructed as a character valued generalization
of the Witten index involving generators of a maximal Cartan subalgebra preserving
one supersymmetry relation. In particular, for a distinguished pair of
supercharges $Q=\overline{Q}_{1 }$ and $Q^{\dag}=-{\overline S}_{1}$,
one has the relation
\begin{equation}
Q Q^{\dag}+Q^{\dag}Q = 2{\mathcal H},\quad Q^2= (Q^{\dag})^2=0,\qquad
\mathcal{H}=H-2\overline{J}_3-3R/2.
\label{susy}\end{equation}
Then, the fermionic generators $Q$ and $Q^{\dag}$ commute with the bosonic operators
$\mathcal{R}=H-R/2$ and $J_3$ and with the maximal torus generators of the flavor group
$F_k$. The latter bosonic operators commute between each other as well.
In lagrangian quantum field theory one works with the fields given by
irreducible representations of the group $G_{\text{full}}$
which are realized as operators acting in the Hilbert space.
All the symmetry generators are then defined as functionals of specific
combinations of these fields.
In this situation the superconformal index is formally defined as the following trace
over the Hilbert space of states \cite{KMMR,Romelsberger1}
\begin{eqnarray}
I(y;p,q) = \text{Tr} \Big( (-1)^{\mathcal F}
p^{\mathcal{R}/2+J_3}q^{\mathcal{R}/2-J_3}
\prod_k y_k^{F_k} e^{-\beta {\mathcal H}}\Big),
\label{Ind}\end{eqnarray}
where $(-1)^{\mathcal F}$ is the $\Z_2$-grading operator for representations of
the $SU(2,2|1)$ supergroup. The variables $p,q,y_k,\beta$ are arbitrary group
parameters whose values are restricted by the condition of convergence of \eqref{Ind}.
Presence of the term $(-1)^{\mathcal F}$ shows that all eigenstates of ${\mathcal H}$
with non-zero eigenvalues drop out from this trace because of the cancellation
of bosonic and fermionic state contributions. It means that the superconformal index
is a weighted sum over BPS states which do not form long  multiplets, $Q\psi=Q^\dag\psi=0.$
Because of that there is no $\beta$-dependence in \eqref{Ind}.

This index was computed heuristically on the basis of physical consideration of
theories on curved background $S^3\times \mathbb{R}$ associated with the radial quantization,
or $S^3\times S^1$ in the Euclidean space.
Space-time symmetry group is reduced and conformal invariance is in general absent
(it emerges in the infrared fixed points). Still, the meaning of operators
entering  \eqref{Ind} as Cartan generators preserving supersymmetry remains intact.

The field theories of interest may contain the vector
superfield  which is always in the adjoint representation of
the gauge group $G$ with the corresponding character $\chi_{adj}(z)$,
and it is invariant with respect to $F$.
They involve also a set of chiral superfields
transforming as certain irreducible representations of the gauge group
with the character $\chi_{R_G,j}(z)$ and of the flavor group $F$
with the characters $\chi_{R_F,j}(y)$ (index $j$ counts such representations).
The antichiral fields are described by conjugated
representations with the characters $\chi_{{\bar R}_G,j}(z)$
and $\chi_{{\bar R}_F,j}(y)$.  The characters depend on the maximal torus variables
$z_a$, $a=1,\ldots,\text{rank}\, G$, and $y_k$, $k=1,\dots,\text{rank}\, F$.

The final result for the index can be represented in the following explicit form:
\begin{equation}
I(y;p,q) \ = \ \int_{G} d \mu(z)\,
\exp \Big( \sum_{n=1}^{\infty}
\frac 1n \ind\big(p^n ,q^n, z^n , y^ n\big ) \Big),
\label{int_G}\end{equation}
where $d \mu(z)$ is the Haar measure for the gauge group $G$ and
\begin{eqnarray}\nonumber  &&
\ind(p,q,z,y) =  \frac{2pq - p - q}{(1-p)(1-q)} \chi_{\adj_G}(z)
 \\ && \makebox[2em]{}
+ \sum_j \frac{(pq)^{R_j/2}\chi_{R_F,j}(y)\chi_{R_G,j}(z) - (pq)^{1-R_j/2}
\chi_{{\bar R}_F,j}(y)\chi_{{\bar R}_G,j}(z)}{(1-p)(1-q)}
\label{1ind} \end{eqnarray}
with some fractional numbers $R_j$ called $R$-charges.
The function \eqref{1ind} is called the one-particle states index
and the integrand of \eqref{int_G} is called the
plethystic exponential. Emergence of the integration over $G$ reflects
the fact that the trace in \eqref{Ind} is taken over the gauge invariant states.

For example, for $G=SU(N)$ one has $z=(z_1,\ldots,z_N)$ with $\prod_{j=1}^Nz_j=1$.
The gauge group measure for functions depending only on $z_j$ has the form
\begin{eqnarray*} &&
\int_{SU(N)} d\mu(z) \ = \   \frac{1}{N!} \int_{\mathbb{T}^{N-1}}
\Delta(z) \Delta(z^{-1}) \prod_{a=1}^{N-1} \frac{dz_a}{2 \pi {i} z_a},
\end{eqnarray*}
where $\Delta(z) \ = \ \prod_{1 \leq a < b \leq N} (z_a-z_b)$.
The fundamental representation character has the form
$\chi_{SU(N),f}(z)=\sum_{k=1}^Nz_k$, and the
adjoint representation character is
$\chi_{SU(N),\adj}(z)= (\sum_{i=1}^N z_i)(\sum_{j=1}^N z_j^{-1})-1$.

Consider the field theory with $(G=SU(2),$ $F=SU(6))$ containing two representations.
The vector superfield transforming as $(\adj, 1)$ with the character
$\chi_{SU(2),\adj}(z)=z^2+z^{-2}+1$. The chiral superfield which is described by the
fundamental representations of both groups $(f, f)$ with the characters
$\chi_{SU(2),f}(z)=z+z^{-1}$ and
$$
\chi_{SU(6),f}(y)=\sum_{k=1}^6y_k,
\quad \chi_{SU(6),\bar f}(y)=\sum_{k=1}^6y_k^{-1},
\quad \prod_{k=1}^6y_k=1.
$$
Let us fix also the chiral field $R$-charge as $R=1/3$.

{\em Exercise:} show that after plugging these data
into the formula \eqref{Ind} and passing from the plethystic exponential to
the infinite product form of the integrand, one obtains precisely the left-hand
side expression for the elliptic beta integral evaluation formula \eqref{elbeta}
after the identification $t_k=(pq)^{1/6}y_k$.

In this picture the unitarity condition
for $SU(6)$ group expressed by the equality $\prod_{k=1}^6y_k=1$ becomes the
balancing condition $\prod_{k=1}^6t_k=pq$ for the integral which is associated
with the hidden ellipticity condition.

Thus, the elliptic beta integral describes
the superconformal index $I_E$ of a particular four dimensional gauge field theory.
Consider now another field theory without gauge group $G=1$ and
containing only one free chiral superfield transforming as the antisymmetric tensor of
the second rank $T_A$ of the same flavor group $F=SU(6)$. The corresponding character is
$$
\chi_{SU(6),T_A}(y)=\sum_{1\leq i<j\leq 6}y_iy_j,
$$
and we fix the $R$-charge for this field as $R=2/3$.

{\em Exercise:} check that substituting these data to the same formula \eqref{Ind}
one comes precisely to the right-hand side expression in \eqref{elbeta}.

So, the result of evaluation of the elliptic beta integral yields the superconformal index
$I_M$ of a completely differently looking
field theory than in the previous case. The two described theories represent the
simplest example of the so-called Seiberg duality \cite{seiberg} which states a
conjectural equivalence of two models in their infrared fixed points.
It is a natural extension of the electromagnetic duality to non-abelian
gauge field theories. Therefore the first described model is called
the ``electric" theory and the second model -- the ``magnetic" one.
The equality of superconformal indices of these two models, $I_E=I_M$, expressed by the evaluation
formula \eqref{elbeta} can be considered as a proof of this duality in the sectors
of BPS states which appear to be identical. The physical phenomenon when the theory
in the ultraviolet regime with nontrivial gauge interaction
becomes in the low energy regime an effective field theory without gauge degrees of freedom
is called the confinement. Thus, the process of computation of the elliptic beta
integral is equivalent to the transition from high to lower energy physics.
From mathematical point of view it describes some group-theoretical duality,
when a particular function on characters  yields the same result
for two different sets of representations of two different groups.

Consider now the full Seiberg electric-magnetic duality \cite{seiberg}.
The electric theory has the gauge group $G=SU(N_c)$ and  the
flavor group $SU(N_f)_l\times SU(N_f)_r\times U(1)_B$
(it enlarges to $SU(2N_f)$ for $N_c=2$). The
representation properties of the fields are described in the table below
(where $\widetilde{N}_c=N_f-N_c$ ):
\begin{center}\begin{tabular}{|c|c|c|c|c|c|}
\hline
& $SU(N_c)$ & $SU(N_f)_l$ & $SU(N_f)_r$ & $U(1)_B$ & $U(1)_R$ \\
\hline
$Q$ & $f$ & $f$ & 1 & 1 & $\widetilde{N}_c/N_f$ \\
$\widetilde{Q}$ & $\overline{f}$ & 1 & $\overline{f}$ & -1 & $\widetilde{N}_c/N_f$ \\
$V$ & ${\rm adj}$  & $1$   &  $1$ &  $0$   &  $1$ \\
\hline
\end{tabular}
\end{center}

The magnetic theory has different gauge group $G=SU(\tilde{N}_c)$ and the same flavor group.
The representation properties of the fields are described in the next table:
\begin{center}
\begin{tabular}{|c|c|c|c|c|c|}
\hline
& $SU(\widetilde{N}_c)$ & $SU(N_f)_l$ & $SU(N_f)_r$ & $U(1)_B$ & $U(1)_R$ \\
\hline
$q$ & $f$ & $\overline{f}$ & 1 & $N_c/\widetilde{N}_c$ & $N_c/N_f$ \\
$\widetilde{q}$ & $\overline{f}$ & 1 & $f$ & $-N_c/\widetilde{N}_c$ & $N_c/N_f$
\\
$M$ & 1 & $f$ & $\overline{f}$ & 0 & $2\widetilde{N}_c/N_f$
\\
$\widetilde V$ & ${\rm adj}$  & $1$   &  $1$ &  $0$   &  $1$ \\
\hline
\end{tabular}
\end{center}

The first columns of these tables contain usual notation for the fields and
last columns contain the abelian group charges -- eigenvalues
of the generators of $U(1)_B$ and $U(1)_R$ groups. The vector superfields
are described in the last rows with all other rows describing some chiral superfields.
According to Seiberg's conjecture, these two $\mathcal{N}=1$ supersymmetric models have
identical physical behaviour at their infrared fixed points
where superconformal symmetry is fully realized. The suggested consistency
checks included the facts that the global anomalies of theories match
('t Hooft anomaly matching conditions) and that the reductions $N_f\to N_f-1$ match for
both theories.
Validity of both criteria can be traced from the equality of
the electric and magnetic theory indices which we describe now.

Superconformal indices for these general theories were constructed in \cite{DO}
(see also \cite{SV}) and we skip the details of their computation.
After passing from maximal torus variables for the flavor group to the
canonical elliptic hypergeometric integral parameters, the electric theory index
takes the form:
\begin{eqnarray}\nonumber
&& I_E = \kappa_{N_c} \int_{\mathbb{T}^{N_c-1}}
\frac{\prod_{i=1}^{N_f} \prod_{j=1}^{N_c} \Gamma(s_i z_j,t_i
z^{-1}_j;p,q)}{\prod_{1 \leq i < j \leq N_c} \Gamma(z_i
z^{-1}_j,z_i^{-1} z_j;p,q)} \prod_{j=1}^{N_c-1} \frac{d z_j}{z_j},
\end{eqnarray}
where $ST = (pq)^{N_f-N_c},\, S=\prod_{i=1}^{N_f} s_i,\, T = \prod_{i=1}^{N_f} t_i,$ and
$$
\prod_{j=1}^{N_c} z_j =1, \qquad
\kappa_{N_c}=\frac{(p;p)_\infty^{N_c-1} (q;q)_\infty^{N_c-1}}{N_c! (2 \pi {i})^{N_c-1}}.
$$
This is a multiple integral for the root system $A_{N_c-1}$, which coincides with
\eqref{AI} for $N_f=N_c+1$ and $n=N_c-1$.

For the magnetic theory one has:
$$
I_M = \kappa_{\tilde{N}_c}\prod_{i,j=1}^{N_f}
\Gamma(s_i t_j;p,q)
  \int_{\mathbb{T}^{\widetilde{N}_c-1}}
\frac{\prod_{i=1}^{N_f} \prod_{j=1}^{\widetilde{N}_c}
\Gamma(S^{\frac{1}{\widetilde{N}_c}} s_i^{-1}
x_j,T^{\frac{1}{\widetilde{N}_c}} t_i^{-1} x_j^{-1};p,q)}
{\prod_{1 \leq i < j \leq \widetilde{N}_c}
\Gamma(x_ix_j^{-1},x_i^{-1}x_j;p,q)}\prod_{j=1}^{\widetilde{N}_c-1}  \frac{dx_j}{x_j},
$$
where $\prod_{j=1}^{\tilde N_c} x_j=1$, $\tilde N_c=N_f-N_c$.

As observed by Dolan and Osborn \cite{DO}, the dual indices coincide $I_E=I_M$,
since the equality of corresponding elliptic hypergeometric integrals
was rigorously established by Rains \cite{rai:trans}
(for some particular values of the parameters it was proven or conjectured
by the author \cite{spi:umn,spi:theta2}). Evidently, this identity is a multivariable
extension of the  second $V$-function transformation law \eqref{E7-2}.

In the case when the electric index is
explicitly computable, i.e. $N_f=N_c+1$, one has the confinement of colored particles
without chiral symmetry breaking.
For $N_f=N_c$ one has the confinement with chiral symmetry breaking which
 is reflected in the appearance of Dirac delta-functions in the description
 of indices \cite{SV3}. In general, equality of dual indices is currently the most
 rigorous mathematical justification of the Seiberg duality conjecture.

Reduction of the number of chiral fields $N_f\to N_f-1$ is reached by the
restriction of the parameters $s_{N_f}t_{N_f}=pq$. In this case $s_{N_f}$ and $t_{N_f}$
disappear from $I_E$ and the rank of the flavor group of electric theory is reduced by one.
In the magnetic theory it is more involved --- a number of poles start to pinch the integration contour
of $I_M$ and the integral starts to diverge, but the vanishing prefactor
$\Gamma(s_{N_f}t_{N_f};p,q)$ makes the product finite with the effective reduction
of ranks of both the magnetic gauge and flavor groups by one, which matches with
the physical picture of \cite{seiberg}.

As to the 't Hooft anomaly matching conditions, they are described by the modified
analogues of the above integrals $I_E$ and $I_M$ \cite{die-spi:unit}.
Define for the electric theory
\begin{eqnarray}\label{el1}
&& I_E^{mod} =  \kappa^{mod}_{N_c} \int_{-\omega_3/2}^{\omega_3/2}
\frac{\prod_{i=1}^{N_f} \prod_{j=1}^{N_c}
{G}(\alpha_i+u_j, \beta_i-u_j;\mathbb{\omega})}
{\prod_{1 \leq i <   j \leq N_c} {G}(u_i-u_j,-u_i+u_j;\mathbb{\omega})}
 \prod_{j=1}^{N_c-1} \frac{du_j}{\omega_3},
\end{eqnarray}
where $\sum_{j=1}^{N_c} u_j =0$,
$$
\kappa_{N_c}^{mod} = \frac{\kappa(\mathbf{\omega})^{N_c-1}}{N_c!}, \qquad
\kappa(\mathbf{\omega}) =-\frac{\omega_3}{\omega_2}
\frac{(p;p)_\infty(q;q)_\infty(r;r)_\infty}
{(\tilde{q};\tilde{q})_\infty}.
$$
and the balancing condition reads
$$
\alpha+\beta = (N_f-N_c) \sum_{k=1}^3\omega_k,
\qquad \alpha =\sum_{i=1}^{N_f} \alpha_i, \quad
\beta =\sum_{i=1}^{N_f} \beta_i.
$$
We denoted the products of modified elliptic gamma functions as
${G}(a,b;\mathbb{\omega}):={G}(a;\mathbb{\omega}){G}(b;\mathbb{\omega})$.
An analogue of $I_M$ has the form
\begin{eqnarray}\label{mag1}
&& I_M^{mod} =  \kappa^{mod}_{\widetilde{N}_c} \prod_{1 \leq i,j \leq N_f}
{G}(\alpha_i+\beta_j;\mathbb{\omega})
\\ \nonumber && \makebox[3em]{} \times \int_{-\omega_3/2}^{\omega_3/2}
\frac{\prod_{i=1}^{N_f} \prod_{j=1}^{\widetilde{N}_c}
{G}(\alpha/\widetilde{N}_c -\alpha_i+v_j,
\beta/\widetilde{N}_c-\beta_i-v_j;\mathbb{\omega})}
{\prod_{1 \leq i <   j \leq \widetilde{N}_c}
{G}(v_i-v_j, -v_i+v_j;\mathbb{\omega})} \prod_{j=1}^{\widetilde{N}_c-1}
\frac{dv_j}{\omega_3},
\end{eqnarray}
where $\widetilde{N}_c=N_f-N_c$ and $\sum_{j=1}^{\widetilde{N}_c}v_j=0$.

{\em Exercise:} show that $I_E^{mod}=I_M^{mod}$ under the conditions
$$
\text{Im}(\alpha_i/\omega_3), \text{Im}((\alpha/\widetilde N_c-\alpha_i)/\omega_3)  <0, \quad
\text{Im}(\beta_i/\omega_3), \text{Im}((\beta/\widetilde N_c-\beta_i)/\omega_3)<0,
$$
when the integration contour in both integrals can be chosen as the straight line segment
connecting $-\omega_3/2$ and $\omega_3/2$.
In a sketchy way, this is reached by substitution of the expression \eqref{modeg'} to
\eqref{el1}, \eqref{mag1} and analysis of the exponential factors
$e^{\varphi_E}$ and $e^{\varphi_M}$ containing sums of $B_{3,3}$-Bernoulli polynomials.
The phase $\varphi_E$ (or $\varphi_M$) looks like a homogeneous cubic polynomial
of the integration variables $u_j$ (or $v_j$) and parameters $\alpha_j,\beta_j,\omega_i$
divided by $\omega_1\omega_2\omega_3$. However, it appears that the integration
variables cancel out in both of them. As a result, $I_E^{mod} =e^{\varphi_E}\tilde I_E$
and $I_M^{mod} =e^{\varphi_M}\tilde I_M$, where the integrals $\tilde I_E$ and $\tilde I_M$ are
obtained from $I_E$ and $I_M$ after the replacements
$s_j\to e^{-2 \pi {i} \alpha_j/\omega_3}$, $t_j\to  e^{-2 \pi {i}
\beta_j/\omega_3}$, $p\to\tilde p,$ and $q\to \tilde r$. Assuming the original
parametrization $s_j=e^{2 \pi {i} \alpha_j/\omega_2}$ and $t_j= e^{2 \pi {i}
\beta_j/\omega_2}$ this boils down to the modular transformation
$(\omega_2,\omega_3)\to (-\omega_3,\omega_2)$ for $I_E$ and $I_M$.
Explicit computation shows that $\varphi_E=\varphi_M$ and this proves the required equality.

For dual field theories the coincidence of $\varphi_E$ and $\varphi_M$ describes the
't Hooft anomaly matching. Namely, each coefficient of their numerator cubic polynomials
corresponds to a particular triangle Feynman diagram involving fermions and particular
gauge or other currents describing global symmetries of the theories.
The above consideration shows that the ratio of kernels of particular elliptic hypergeometric integrals
corresponding to electric and magnetic superconformal indices has a particular behaviour
from the viewpoint of $\textrm{SL}(3,\mathbb{Z})$-group. One can formalize this statement
in a general setting by taking the following parametrization for such a ratio
\begin{equation}
\Delta(x_1,\dots,x_n;p,q) = (p;p)_\infty^{r_{-}}(q;q)_\infty^{r_{-}}\prod_{a=1}^K
\Gamma\Bigl((pq)^{\frac{R_a}{2}}x_1^{m_1^{(a)}}x_2^{m_2^{(a)}} \dots
x_n^{m_n^{(a)}};p,q\Bigr)^{\epsilon_a},
\label{gen-term}\end{equation}
where $K$ is the total number of independent elliptic gamma functions appearing
in this ratio in the integer power $\epsilon_a$ with its own $R$-charge $R_a$
and $m_j^{(a)}$ -- integer
powers of $n$ independent group parameters $x_j$ (playing the role of fugacities $y_j$ in the
original definition of the superconformal indices). For the Seiberg duality the integer number
$r_{-}$ is equal to the difference between ranks of the electric and magnetic gauge groups.

Using the parametrization $x_j=e^{2\pi{i}u_j/\omega_2}$ one can define a
modified elliptic gamma function analogue of \eqref{gen-term}
\begin{equation}
\Delta^{mod}(u_1,\dots,u_n;\mathbf{\omega}) =  \kappa(\mathbf{\omega})^{r_{-}}
\prod_{a=1}^K {G}\Bigl(R_a\sum_{k=1}^3\frac{\omega_k}{2}
+\sum_{j=1}^nu_jm_j^{(a)};\mathbf{\omega}\Bigr)^{\epsilon_a}.
\label{gen-modterm}\end{equation}

Now one demands validity of an $\textrm{SL}(3,\Z)$-modular transformation relation between
functions \eqref{gen-term} and \eqref{gen-modterm}
\begin{equation}
\Delta^{mod}(u_1,\dots,u_n;\mathbf{\omega}) = \Delta(e^{-2\pi{i}u_1/\omega_3},
\ldots,e^{-2\pi{i}u_n/\omega_3};\tilde p,\tilde r).
\label{modularity}\end{equation}
There are six independent in form functional combinations of $u_j$ and $\omega_i$ in the
sum of $B_{3,3}$-polynomials, appearing after substitution of relation \eqref{modeg'}
in \eqref{modularity}, and additional terms generated by the Dedekind function
modular transformation. The coefficients in front of them should vanish, which
yields the following set of equations
\begin{eqnarray}\label{mcon1}
&& \sum_{a=1}^K \epsilon_a m_i^{(a)} m_j^{(a)} m_k^{(a)} = 0,\\
&& \sum_{a=1}^K \epsilon_a m_i^{(a)} m_j^{(a)}(R_a-1)= 0,
\label{mcon2} \\
&& \sum_{a=1}^K \epsilon_a m_i^{(a)}(R_a-1)^2 =0,
 \label{mcon3}
\\
&& \sum_{a=1}^K \epsilon_a m_i^{(a)} =0,
 \label{mcon4}
\\
&& \sum_{a=1}^K \epsilon_a (R_a-1)^3+r_{-}=0,
 \label{mcon5}
\\
&& \sum_{a=1}^K \epsilon_a (R_a-1)+r_{-} =0.
 \label{mcon6}
\end{eqnarray}
Assuming rationality of $R_a$ we come to a system of Diophantine equations
which were not systematically investigated yet from mathematical point of view, although all
known physical dualities satisfy them as the 't Hooft anomaly matching conditions.
We do not describe the physical meaning of each type of the above equations
referring for details to \cite{SV2}.
We only mention that in the context of superconformal indices the
combinations of integration variables  entering the Bernoulli polynomials must cancel
independently for electric and magnetic indices to be able to pull exponentials $e^{\varphi_{E,M}}$
out of the integrals.

{\em Exercise:} suppose that \eqref{gen-term} is a kernel of an elliptic hypergeometric integral
with $x_1,\ldots,x_r$ being the integration variables, i.e. that it satisfies a set of $r$ $q$-difference
equations in these variables with $p$-elliptic function coefficients.  Show that this requirement
is equivalent to equations \eqref{mcon1} and \eqref{mcon2} with $1\leq i,j\leq r$
together with an extra requirement $\sum_{a=1}^K \epsilon_a m_i^{(a)} m_j^{(a)}\in 2\mathbb{Z}$.

In all known dual theories the latter extra evenness condition is automatically satisfied, though it is
not clear whether it follows from general equations \eqref{mcon1}-\eqref{mcon6}. Condition \eqref{mcon1} for
all $1\leq i,j,k\leq r$ physically corresponds  to the demand of absence of the
gauge anomalies, which is needed for the consistency of field
theories whose indices are described by the corresponding integrals.
As we see, it follows from the original definition of the elliptic hypergeometric integrals
 \eqref{EHImult} and its multivariable extension, which thus gets an interesting physical interpretation.

As a summary of connections with the superconformal indices, we mention that very many
identities for elliptic hypergeometric integrals were found following the physical duality
conjectures, and they still require rigorous proofs, see, e.g. \cite{SV}.
Vice versa, there is a good number of new physical dualities conjectured from proven
integral identities. There are also applications of superconformal indices
to topological field theories, description of lower and higher dimensional field theories,
and some other constructions of mathematical physics \cite{Kim,RR2016}.

\smallskip
{\bf Acknowledgements.}
The author is indebted to E. M. Rains and S. O. Warnaar for helpful discussions.
This work is supported in part by the Laboratory of
Mirror Symmetry NRU HSE, RF government grant, ag. no. 14.641.31.0001.

\end{document}